\documentclass[11pt, a4paper,leqno]{amsart}
\usepackage{amsmath,amsthm,amscd,amssymb,amsfonts, amsbsy}
\usepackage{latexsym}
\usepackage{txfonts}
\usepackage{exscale}
\usepackage{mathtools}
\usepackage{mathrsfs}
\usepackage[titletoc]{appendix}

\usepackage{amssymb}

\usepackage[utf8]{inputenc}

\usepackage[colorlinks,linkcolor=blue,citecolor=blue,pagebackref,hypertexnames=false, breaklinks]{hyperref}
\usepackage[usenames,dvipsnames]{color}

\calclayout
\allowdisplaybreaks

\def \Int{\,\Int\,}

    \def\XXint#1#2#3{{\setbox0=\hbox{$#1{#2#3}{\int}$}
    \vcenter{\hbox{$#2#3$}}\kern-.5\wd0}}

    \def\XXint#1#2#3{{\setbox0=\hbox{$#1{#2#3}{\int}$}
    \vcenter{\hbox{$#2#3$}}\kern-.5\wd0}}

\def \C{ \mathbb{C} }

\def \R{ \mathbb{R} }
\def \Z{ \mathbb{Z} }

\renewcommand{\chi}{\mathbf{1}}

\theoremstyle{plain}
\newtheorem{theorem}{Theorem}[section]
\newtheorem{lemma}{Lemma}[section]
\newtheorem{corollary}{Corollary}[section]
\newtheorem{proposition}{Proposition}[section]

\theoremstyle{definition}
\newtheorem{definition}{Definition}[section]

\newtheorem{example}{Example}[section]

\theoremstyle{remark}
\newtheorem{remark}{Remark}[section]
\begin{document}
\allowdisplaybreaks
\author{Alberto Lastra}
\address{Alberto Lastra
\\
Universidad de Alcalá 
\\
Departamento de Física y Matemáticas
\\
Campus universitario 
\\E-28805 Alcalá de Henares (Madrid), Spain
}
 \email{alberto.lastra@uah.es}
\author{Cruz Prisuelos-Arribas}
\address{Cruz Prisuelos-Arribas
\\
Universidad de Alcalá 
\\
Departamento de Física y Matemáticas
\\
Campus universitario 
\\E-28805 Alcalá de Henares (Madrid), Spain
} \email{cruz.prisuelos@uah.es}

\author{V\'ictor Soto-Larrosa}
\address{V\'ictor Soto-Larrosas
\\
Universidad de Alcalá 
\\
Departamento de Física y Matemáticas
\\
Campus universitario 
\\E-28805 Alcalá de Henares (Madrid), Spain
} \email{v.soto@uah.es}

\title[]{On linear systems of moment differential equations with singularities of the first kind}


\date{\today}

\begin{abstract}
The solution to systems of moment differential equations of the form $z\partial_my=(zA+B)y$ are provided, for a matrix $B$ with general good spectrum. Existence and convergence of Floquet-type solutions is studied. A generalized definition of $z^B$ is given, as a tool to solve the main problem whenever $A\equiv 0$. The theory is illustrated with examples which are important in applications. 
\end{abstract}

\maketitle

\tableofcontents

\section{Introduction}\label{secintro}
In this paper we consider  systems of moment differential equations of the form
\begin{align}\label{problem}
 z\partial_my(z)=(zA+B)y(z)
\end{align}
where $A, B\in \mathbb{C}^{n\times n}$ are  constant matrices and $y(z)=(y_1(z),\ldots,y_n(z))$ represents a vector of unknown functions, for some positive integer $n\geq 1$. The symbol $\partial_my(z)=(\partial_m(y_1),\ldots,\partial_m(y_n))$ stands for the moment derivative of $y(z)$. 

The study of moment derivatives and moment differential equations has attracted the attention of many researchers during the last decade. Its versatility relies on its different realizations, depending on the choice of the moment sequence $m$. In this sense, moment derivative $\partial_m$is shown as the usual derivative, a fractional derivative, Jackson $a-$derivative, etc. depending on the choice of $m$. On the one hand, valid results attained in the framework of moment differential equations/systems remain valid for all the different realizations of the moment derivative at once. On the other hand, such a level of generality causes substantial difficulties when handling it. For example, the operator $\partial_m$ lacks of a Leibniz-like rule.

The first time that the moment differential operator was considered was in the work~\cite{balseryoshino} of 2010, by W. Balser and M. Yoshino, associated to techniques of kernel functions for summability of Gevrey type. Afterwards, kernel functions for generalized summability were developed in~\cite{sanz}, and its associated generalized summability in~\cite{lastramaleksanz}. In principal, moment differentiation is initially formally defined on formal power series or holomorphic functions near the origin  by means of their Taylor expansion. However, moment differentiation has also been extended to functions which are the generalized sum (resp. generalized multisum) of formal solutions to moment differential equations. We refer to~\cite{lastramichaliksuwinska21} (resp.~\cite{lastramichaliksuwinska23}). It is also worth mentioning the work~\cite{su} providing Gevrey estimates on these type of functional equations, and the study of Stokes phenomenon in~\cite{mitk}. See the references in the mentioned studies for further advances on the topic.

Asymptotic results related to generalized summability of formal power series have also been a branch of study due to the appearance of a great variety of phenomena which goes beyong Gevrey classes. For example, the sequence $1^{+}-$Gevrey sequence is given by $m(0)=m(1)=1$ and $m(n)=(n/\log(n))^{n}$ for every $n\ge2$. This sequence is of great importance in asymptotic theory due to its  appearance as natural bounds when relating the analytic and formal solution to functional equations, mainly difference equations, see~\cite{im,lastramalek19,ma14}. Further advances on this topic study the regularity of solutions to moment (partial) differential equations in~\cite{resu}.

In the last years, systems of moment differential equations have also been studied, following the knowledge of systems of differential equations. In~\cite{lastra}, the first author describes the general solution of a linear system of moment differential equations 
$$\partial_my=Ay,$$ 
with $A\in\C^{n\times n}$. In this work, the growth rate of the entire solutions at infinity is also described from different points of view. Lately, in~\cite{lastraprisuelos}, a compact description of the solutions to such systems is provided by means of a generalized definition of the exponential of a matrix. In this work, we give a step forward in this direction by considering systems of moment differential equations whose matrix of coefficients is of the form $\frac{1}{z}B+A$, for some $A,B\in\C^{n\times n}$, adding the possibility of a simple pole at the origin with respect to the previous situation. This is the most natural advance when following the classical study of systems of linear systems of differential equations in the complex domain: see for example Chapter II,~\cite{wa} or Chapter 2 in~\cite{balser}. We also refer to~\cite{deligne,yoshida} about the classical theory of singularities of first kind.

The main advances in the present work apply, but are not restricted, to moment derivatives associated to sequences of moments determined by kernels for generalized summability (see Section~\ref{sec2}). Most of the results provided are also valid when considering sequences of positive real numbers which turn out to be the restriction of a function $m:\{z\in\C:\hbox{Re}(z)\ge0\}\to\C$, to the set of nonnegative integers. Indeed, we focus our attention to particular sequences of moments which are determined by the mentioned theory, such as $(\Gamma(1+\frac{\alpha}{p}))_{p\ge0}$ for some $\alpha>0$ (related to Caputo fractional derivative as the realization of the moment derivative), but also to sequences which fall appart from the previous theory such as the sequence of $q-$factorials, $([p]_{q}^{!})_{p\ge0}$ (for which the moment derivative coincides with Jackson's $q$-derivative).

The study of systems of the form (\ref{problem}) makes it necessary to introduce novel elements and tools describing their solution. In this regard, we extend the definition of $\partial_m$ to formal power series of the form $z^{\nu}\C[[z]]$, for some $\nu\in\C$ with $\hbox{Re}(\nu)\ge1$ (see Definition~\ref{defi385}).
In Section~\ref{sec3}, we provide formal and analytic solutions to (\ref{problem}) under Assumption (H1), i.e. the matrix $B$ having a ``good spectrum'' in the terminology of the classical theory (Theorem~\ref{teopral1}). As a novel tool describing the solution to (\ref{problem}), we define the matrix $z_m^{B}$, as a generalization of the matrix $z^{B}$, solution to $z\partial_my=By$ in the case that $B$ is a diagonalizable matrix. The results obtained are refined substituting the asumptions for existence and convergence of formal solutions in the form of a Floquet series in Section~\ref{sec35} and describe some important consequences of the theory to applications to systems of fractional differential equations, and to systems of $q$-difference equations (Section~\ref{secapap}). 

On the way, some results on the approach to the problem (\ref{problem}) are provided: the simplification of the system and the relation between the initial and modified solutions (Lemma~\ref{lematec}) together with their convergence; msimplification of $B$ to its Jordan canonical form (Proposition~\ref{prop1}). These results are applied to show the construction of two solutions in the particular case of planar systems (Section~\ref{secdiag} and Section~\ref{secnodiag}). Finally, the absence of two Flowuet solutions in the case of $B$ consisting of a unique Jordan block motivates the introduction of an additional hypothesis on the existence of elements satisfying a moment differential equation in order to obtain a second solution to the planar system (see (\ref{e1141})). This fact motivates the statement of Assumption (H3) in the higher dimensional case to provide the solution to $z\partial_m y=By$ in the form of $z_m^{B}$, whenever $B$ is not diagonalizable (recall that $z_m^B$ was already defined in the diagonalizable case). Theorem~\ref{teopral2} describes the construction of $z_m^B$ as a solution to the problem. Finally, we justify that Assumption (H3) is feasible in concrete situations.

\vspace{0.3cm}

\textbf{Notation:}

We write $\C[[z]]$ for the set of formal power series with complex coefficients in the variable $z$.

For any fixed open domain $U\subseteq\C$, we write $\mathcal{O}(U)$ for the set of holomorphic functions defined in $U$.

Let $d\in\R$ and $0<\theta<2\pi$. We write $S_{d}(\theta)$ for the infinite sector $S_d(\theta)=\{z\in\C:|\hbox{arg}(z)-d|<\theta/2\}$, and for given $z_0\in\C$ and positive $r$, we write $D(z_0,r)$ for the open disc $\{z\in\C:|z-z_0|<r\}$.

Given $A\in\C^{n\times n}$, $\hbox{spec}(A)$ stands for the spectrum of $A$. 

\section{Moment sequences and moment derivatives}\label{sec2}

In this section, we recall the main facts on the operator of moment derivation, generalizing the classical derivative. This formal operator is constructed by means of a so-called sequence of moments, say $m$. Such sequence is, in principal, a sequence of positive real numbers. However, in most applications $m$ often emerges as the sequence of moments associated to some measure. We will mainly focus on families of moments associated to a pair of kernel functions for generalized summability.

Let us consider a sequence of positive real numbers $m=(m_p)_{p\ge0}$, and define the formal operator $\partial_m:\C[[z]]\to\C[[z]]$ by
\begin{equation}\label{e297}
\partial_m\left(\sum_{p\ge0}\frac{a_p}{m_p}z^p\right)=\sum_{p\ge0}\frac{a_{p+1}}{m_p}z^p.
\end{equation}
Several examples are at hand:
\begin{itemize}
\item The most outstanding example is the usual derivative, which appears when considering the sequence $m=(p!)_{p\ge0}$.
\item Given $k>0$, the moment differential operator determined by the sequence $m_{k}=(\Gamma(1+\frac{p}{k}))_{p\ge0}$ is linked to Caputo fractional derivative by the relation
$$\partial_{m_k}f(z^{1/k})=\prescript{C}{}{D_z}^{1/k}(f(z^{1/k})),$$
when applied to a formal power series. We refer to~\cite{michalik12}, Def. 5 and Rem. 1 for further details.
\item Let $q\in(0,\infty)\setminus\{1\}$. The moment derivative associated with the sequence of $q-$factorials is known as the $q$-Gevrey sequence $m_q=([p]_q!)_{p\ge0}$, with $[\ell]_q!:=[1]_q[2]_q\cdots[\ell]_{q}$ for every positive integer $\ell$, and $[0]!_{q}=1$. Here, $[j]_q=\sum_{n=0}^{j-1}q^{n}$ for $j\ge 1$. This moment derivative coincides from a formal point of view with the $q-$derivative $D_q$ defined by
$$D_qf(z)=\frac{f(qz)-f(z)}{qz-z}.$$
\end{itemize}

The concept of moment derivation goes back to the work~\cite{balseryoshino}, by W. Balser and M. Yoshino. In that work, the sequence of moments is constructed by means of a pair of kernel functions (see~\cite{balser}, Section 6.5), remaining close to the classical methods applied in the theory of summability of formal solutions to functional equations. Indeed, such methods are equivalent to the classical ones, as there exist positive constants $A_jC_j$, $j\in\{1,2\}$, such that
$$C_1A_1^{p}\Gamma(1+\frac{p}{k})\le m_p\le C_2A_2^{p}\Gamma(1+\frac{p}{k}),$$
for every $p\ge0$. The previous methods were extended by J. Sanz and coauthors in~\cite{sanz,lastramaleksanz,jimenezkamimotolastrasanz} to a more general framework, which involves the existence of a pair of kernel functions for generalized summability. In the following paragraphs, we recall the notion of pair of kernel functions for generalized summability, which determines the key point in such generalized methods. 

\begin{definition}[\cite{thilliez}, Section 1.1]
Let $\mathbb{M}=(M_p)_{p\ge0}$ be a sequence of positive real numbers, normalized to $M_0=1$, and satisfying the following properties:
\begin{itemize}
\item[(lc)] $M_p^2\le M_{p-1}M_{p+1}$ for $p\ge1$ ($\mathbb{M}$ is logarithmically convex)
\item[(mg)] $M_{p+1}\le A^{p+q}M_pM_q$ for some $A>0$, valid for all $p,q\ge0$ ($\mathbb{M}$ is a sequence of moderate growth)
\item[(snq)] There exists $B>0$ with
$$\sum_{q\ge p}\frac{M_q}{(q+1)M_{q+1}}\le B\frac{M_p}{M_{p+1}},$$
for every $p\ge 0$.
\end{itemize}
Then, $\mathbb{M}$ is a strongly regular sequence.
\end{definition}
 
The sequence $m_{k}$ above turns out to be an example of a strongly regular sequence. In particular, the sequence of factorials is a strongly regular sequence. Gevrey sequence of order $1/k$, defined by $(p!^{1/k})_{p\ge0}$, is also a strongly regular sequence. The sequence of $q-$factorials is not a strongly regular sequence, as it does not satisfy (mg) if $q>1$, and it does not safisfy (snq) if $0<q<1$.

For any strongly regular sequence $\mathbb{M}=(M_p)_{p\ge0}$, we consider the function
$$M(t)=\sup_{p\ge0}\log\left(\frac{t^p}{M_p}\right)$$ 
for every $t>0$, and $M(0)=0$. It turns out that $M$ defines a positive function in $(0,\infty)$ for which
$$\rho(M):=\lim\sup_{t\to\infty}\max\left\{0,\frac{\log(M(t))}{\log(t)}\right\}$$
determines a positive real number. Associated to a strongly regular sequence, one has at hand the notion of pair of kernel functions for generalized summability.

\begin{definition}\label{defi332}
Let $\mathbb{M}$ be a strongly regular sequence with $\rho(M)>1/2$. We say that $\mathbb{M}$ admits a pair of kernel functions for $\mathbb{M}$-summability if there exists a pair $(e,E)$ such that:
\begin{itemize}
\item[(i)] $e\in\mathcal{O}(S_0(\pi/\rho(M)))$ which is locally uniformly integrable at the origin, i.e. there exists $t_0>0$ and for all $z_0\in S_0(\pi/\rho(M))$ there exists $r_0>0$ with $D(z_0,r_0)\subseteq S_0(\pi/\rho(M))$ such that
$$\int_0^{t_0}\sup_{z\in D(z_0,r_0)}|e(t/z)|\frac{dt}{t}<\infty.$$
The function $e$ is subject to an exponential decay at infinity: for every $0<\Delta<\pi/\rho(M)$ there exist $C_1,C_2>0$ with
$$|e(z)|\le C_1\exp(-M(|z|/C_2)),\qquad z\in S_0(\Delta).$$
\item[(ii)] $e$ restricted to the positive real line defines a real positive function.
\item[(iii)] The function $E$ is entire, and of exponential growth at infinity with growth governed in terms of 
$$|E(z)|\le C_3\exp(M(|z|/C_4)),\qquad z\in\C$$
for some $C_3,C_4>0$. In addition to this, there exists $\beta>0$ such that for all $0<\theta<\pi(2-1/\rho(M))$ and $R>0$ there exists $C_5>0$ with
$$|E(z)|\le \frac{C_5}{|z|^{\beta}},\qquad z\in S_{\pi}(\theta)\cap(\C\setminus D(0,R)).$$ 
\end{itemize}
Both functions are related via the moment function associated with $e$, defined by
$$m_e(z)=\int_0^{\infty}t^{z-1}e(t)dt,$$
for every $z\in\C$ with $\hbox{Re}(z)\ge0$. 

It is worth mentioning that the previous assumptions on $e$ guarantee that the moment function is well-defined and continuous in $\{z\in\C:\hbox{Re}(z)\ge0\}$, holomorphic in $\{z\in\C:\hbox{Re}(z)>0\}$, with $m_e(x)>0$ for all $x\ge0$. 

The sequence $(m_e(p))_{p\ge0}$ defined by
\begin{equation}\label{e348}
m_e(p)=\int_0^{\infty}t^{p-1}e(t)dt,\qquad p\ge0
\end{equation}
 is known as the sequence of moments associated with the pair of kernel functions.

On the other hand, the entire function $E$ is represented by the following Taylor expansion at the origin:
\begin{equation}\label{e348b}
E(z)=\sum_{p\ge0}\frac{1}{m_e(p)}z^p,\qquad z\in\C.
\end{equation}
\end{definition}

\begin{remark} A convenient ramification of the variable allows to adapt the previous definition to a strongly regular sequence with $0<\rho(M)\le1/2$. See~\cite{lastramaleksanz}, Remark 3.5 (iii). In our settings, the representation of the entire function $E$ in (\ref{e348b}) remains valid, and therefore no distinction is made.
\end{remark}

It is important to point out that not every strongly regular sequence $\mathbb{M}$ admits a pair of kernel functions. However, this is guaranteed in the case that $\mathbb{M}$ admits a nonzero proximate order in the sense of Lindel\"of (see~\cite{jss}). Strongly regular sequences which usually appear in applications do satisfy this requirement. For example, given $0<\alpha<2$, Gevrey sequence of order $\alpha$, $\mathbb{M}_{\alpha}=(p!^\alpha)_{p\ge0}$, admits the pair of kernel functions $(e_{\alpha},E_{\alpha})$, with 
$$e_\alpha(z)=\frac{1}{\alpha}z^{1/\alpha}\exp(-z^{1/\alpha}),\quad z\in S_{0}(\alpha\pi),$$
and its associated sequence of moments $(m_{e_{\alpha}}(p))_{p\ge0}$ is given by $m_{e_{\alpha}}(p)=\Gamma(1+\alpha p)$, leading to
$$E_{\alpha}(z)=\sum_{p\ge0}\frac{1}{\Gamma(1+\alpha p)}z^p,\qquad z\in\C.$$

In particular, moment derivatives in the form (\ref{e297}) can be applied when considering a sequence $m$ being the sequence of moments associated to a pair of kernel functions for generalized summability, as showed in Definition~\ref{defi332}. This is the cornerstone of the theory of generalized summability of formal solutions to moment differential equations, in which the action of the moment differential operator acts on some functional equation under study. As a consequence, any result proved for a moment differential equation is valid for a plethora of situations, including fractional differential equations and usual differential equations.  

The definition of moment differentiation, initially of formal nature, can be naturally extended to holomorphic functions defined on some neighborhood of the origin, in particular to entire functions, by identification of the functions with their Taylor expansion at the origin. Additionally, the definition of moment differentiation has already been extended to the set of holomorphic functions defined on sectors of the complex domain with vertex at the origin which turn out to be the generalized sum~\cite{lastramichaliksuwinska21}, or the generalized multisum~\cite{lastramichaliksuwinska23} of some formal power series. We refer to the previous references for more details on the previous facts.

However, there are other important cases of moment derivatives, which are not associated to a strongly regular sequence, admitting pairs of functions with analogous properties as the pair of kernel functions for generalized summability. The most important example is that of the $q$-Gevrey sequence, for given $q>1$. This is not a strongly regular sequence as its growth rate surpasses that of a sequence satisfying (mg) condition. For such sequence, the function $\sum_{p\ge0}\frac{1}{[p]!_q}z^p$ is known in the literature as one of the $q$-exponential functions. The $q-$Gevrey sequence $([p]_q!)_{p\ge0}$ is intimately related to $(q^{p(p-1)/2})_{p\ge0}$, due to
$$\lim_{p\to\infty}\frac{[p]_q!}{q^{\frac{p(p-1)}{2}}}\left(1-\frac{1}{q}\right)^{n}=c(q),$$
for some constant $c(q)>0$ only depending on $q$. Therefore, a formal power series $\sum_{p\ge0}\frac{a_p}{[p]_q!}z^{p}\in\C[[z]]$ has a positive radius of convergence if and only if the same holds for the formal power series $\sum_{p\ge0}\frac{a_p}{q^{\frac{p(p-1)}{2}}}z^{p}\in\C[[z]]$. The sequence $(q^{p(p-1)/2})_{p\ge0}$ is the sequence of moments (in accordance to~(\ref{e348})) associated to the kernel function $e(z)=\sqrt{2\pi\ln(q)}\exp\left(\frac{\log^{2}(\sqrt{q}z)}{2\ln(q)}\right)$, and also $e(z)=\ln(q)\Theta_{1/q}(z)$, with $\Theta_{1/q}$ being Jacobi Theta function. See~\cite{carrillolastra,ramis,viziozhang} for more details. Another example of explicit construction of kernel functions associated to a sequence which is nos strongly regular can be found in~\cite{jls}.

We extend the definition of a moment derivative to formal power series in the monomials $\{z^{p+\nu}:p\ge0\}$, for some fixed $\nu\in\C$ with $\hbox{Re}(\nu)\ge1$.  This object will allow us to solve certain meromorphic linear systems of moment differential equations. In principal, the previous monomials will be considered from a formal point of view, leading to symbolic solutions, but also associated to the function $z^c=\exp(c\log(z))$ (analytic on some domain when fixing some branch of the complex logarithm, or as a multi-valued function.

\begin{definition}\label{defi385}
Let $\nu\in\C$ with $\hbox{Re}(\nu)\ge1$. Let $m$ be the moment function associated to some pair of kernel functions $(e,E)$, as described in Definition~\ref{defi332}. We extend the definition of $\partial_m$ to the set 
$$z^{\nu}\C[[z]]=\{\sum_{p\ge0}a_pz^{\nu+p}: a_p\in\C, p\ge0\}$$ 
by
\begin{equation}\label{e390}
\partial_m\left(\sum_{p\ge0}a_pz^{p+\nu}\right)=\sum_{p\ge0}a_p\frac{m(p+\nu)}{m(p+\nu-1)}z^{p+\nu-1}.
\end{equation}
\end{definition}

Observe the previous definition makes sense due to $\hbox{Re}(\nu -1)=\hbox{Re}(\nu)-1\ge0$, where $m$ is well-defined. It is also worth mentioning that Definition~\ref{defi385} can be naturally extended to any sequence $m$ of positive real numbers.

It is worth mentioning that the previous definition is linked to the process of searching for solutions of functional equations, known as Frobenius method. The following result is a direct consequence of the previous definition.

\begin{proposition}
In the situation of Definition~\ref{defi385}, the following statements hold:
\begin{itemize}
\item $\partial_m(\hat{f})\in z^{\nu-1}\C[[z]]$ for every $\hat{f}\in z^{\nu}\C[[z]]$.
\item Let $\nu$ be a positive integer. Then, Definition~\ref{defi385} of $\partial_m$ coincides with that of (\ref{e297}).
\item Let $m=(p!)_{p\ge0}$. Then, Definition~\ref{defi385} of $\partial_m$ coincides with the usual derivative of a generalized complex power function.
\end{itemize}
\end{proposition}
\begin{proof}
The first statement is a direct consequence of the definition. 

Let $\nu$ be a positive integer. We observe from the definition of moment derivative that
$$\partial_m\left(\sum_{p\ge0}a_pz^{p+\nu}\right)=\partial_m\left(\sum_{p\ge\nu}a_{p-\nu}z^{p}\right)=\sum_{p\ge\nu}\frac{a_{p-\nu}m(p)}{m(p-1)}z^{p-1},$$
which coincides with (\ref{e390}).

Finally, we observe that 
\begin{multline*}
\left(\sum_{p\ge0}a_pz^{p+\nu}\right)'=\left(z^{\nu}\sum_{p\ge0}a_pz^{p}\right)'=\nu z^\nu\sum_{p\ge0}a_pz^{p-1}+z^{\nu}\sum_{p\ge1}a_ppz^{p-1}\\
=\nu a_0z^{\nu-1}+\sum_{p\ge1}(\nu+p)a_pz^{\nu+p-1}=\sum_{p\ge0}(\nu+p)a_pz^{\nu+p-1},
\end{multline*}
which concludes the proof.
\end{proof}

%%%%%%%%%%%%%%%%%%%%%%%%%%%%%%%%%%%%%%%%%%%%%%%%%%%%%%%%%%%%%%%%%%%%5
\section{Confluent hypergeometric systems: formal and convergent solutions}\label{sec3}

Let $n\ge1$ be an integer, and fix a sequence of positive real numbers, $m=(m(p))_{p\ge0}$. We consider the linear system of moment differential equations 
\begin{align}\label{eqhypsys}
z\partial_my=(zA+B)y,
\end{align}
for $A,B\in\C^{n\times n}$. 

From now on, we assume that the sequence $m$ is the restriction to the nonnegative integers of some function defined on $\{z\in\C:\hbox{Re}(z)\ge0\}$. We preserve the notation $m=m(z)$ for such function. Observe that this assumption holds whenever one departs from the sequence of moments associated to a pair of kernel functions for generalized summability, or the sequence of moments associated to the $q$-Gevrey sequence, as described in the previous section.

First, we provide formal solutions to (\ref{eqhypsys}). We adopt the following hypothesis:
\begin{itemize}
\item[(H1)] There exists $\mu \in \mathbb{C}$ with $Re(\mu)\geq 1$ such that $\frac{m(\mu)}{m(\mu-1)}\in\hbox{spec}(B)$ and $\frac{m(p+\mu)}{m(p+\mu-1)}$ does not belong to the spectrum of $B$, for all $p\geq 1$.
\end{itemize}

\begin{lemma}\label{lema1}
Let (\ref{eqhypsys}) be a linear system of moment differential equations such that hypothesis (H1) holds. Assume that $s_0$ is an eigenvector associated to the eigenvalue $\frac{m(\mu)}{m(\mu-1)}$ of $B$. Then, the formal power series 
\begin{equation}\label{e416}
y(z)=\sum_{p=0}^{\infty}s_pz^{p+\mu},
\end{equation}
where
\begin{align}\label{solution hypergeometric system}
s_p=\left(\frac{m(p+\mu)}{m(p+\mu-1)}I-B\right)^{-1}As_{p-1},\quad p\geq 1.
\end{align}
is a formal solution to (\ref{eqhypsys}).
\end{lemma}
\begin{proof}
After plugging a formal power series of the form (\ref{e416}) into (\ref{eqhypsys}), and taking into account (H1), one arrives at the recurrence formula 
$$\frac{m(p+\mu)}{m(p+\mu-1)}s_p=As_{p-1}+Bs_p,\qquad p\geq 1,$$
departing from the identity $Bs_0=\frac{m(\mu)}{m(\mu-1)}s_0$, the second part of hypothesis (H1) allows us to conclude the result. Indeed, one has
\begin{multline*}
z\partial_my(z)=z\sum_{p=0}^{\infty}s_p\frac{m(p+\mu)}{m(p+\mu-1)}z^{p+\mu-1}=\sum_{p=0}^{\infty}s_p\frac{m(p+\mu)}{m(p+\mu-1)}z^{p+\mu}
\\
=Bs_0z^{\mu}+\sum_{p=1}^{\infty}s_p\frac{m(p+\mu)}{m(p+\mu-1)}z^{p+\mu}=Bs_0z^{\mu}+\sum_{p=1}^{\infty}(As_{p-1}+Bs_p)z^{p+\mu}
\\
=\sum_{p=1}^{\infty}As_{p-1}z^{p+\mu}+\sum_{p=0}^{\infty}Bs_pz^{p+\mu}=z\sum_{p=0}^{\infty}As_{p}z^{p+\mu}+\sum_{p=0}^{\infty}Bs_pz^{p+\mu}=(zA+B)y(z).
\end{multline*}
\end{proof}

\begin{remark} Observe that the system (\ref{eqhypsys}) admits a formal solution in the form of a Floquet solution (\ref{e416}) if and only if hypothesis (H1) holds.
\end{remark}

We fix any norm on $\C^{n}$, say $\left\|\cdot\right\|$, and consider (maintaining its notation) its matrix induced norm which turns out to be submultiplicative, i.e. 
$$\left\|A B\right\|\le \left\|A\right\|\,\, \left\|B\right\|,$$
for every $A,B\in\C^{n\times n}$.

The following additional assumption provides convergence of formal solutions in the form of a Floquet series. 

\begin{itemize}
\item[(H2)]  There exists $C>0$ such that 
$$\sup_{p\ge1}\left\|\left(\frac{m(p+\mu)}{m(p+\mu-1)}I- B\right)^{-1}\right\|\le C.$$
\end{itemize}

%\begin{lemma}
%For every regular matrix $A\in\C^{n\times n}$ it holds that
%$$1\le \left\|A\right\|\left\|A^{-1}\right\|.$$
%\end{lemma}
%\begin{proof}
%A direct application of the definition of a matrix norm yields
%$$1=\sup_{v\in\C^{n},\left\|v\right\|=1}\frac{\left\|I_nv\right\|}{\left\|v\right\|}=\sup_{v\in\C^{n},\left\|v\right\|=1}\frac{\left\|AA^{-1}_nv\right\|}{\left\|v\right\|}\le \left\|A\right\|\sup_{v\in\C^{n},\left\|v\right\|=1}\frac{\left\|A^{-1}_nv\right\|}{\left\|v\right\|}=\left\|A\right\|\left\|A^{-1}\right\|.$$
%\end{proof}

\begin{corollary}\label{coro0}
Under the hypotheses of Lemma~\ref{lema1}, assume $(H2)$ holds. Then, the system (\ref{eqhypsys}) admits a unique convergent solution in the form (\ref{e416}). 
\end{corollary}
\begin{proof}
Convergence of the solution near the origin is determined by the convergence of  the series $\sum_{p=0}^{\infty}s_pz^p$.
Taking into account the nature of the norm $\left\|\cdot\right\|$, we observe from %Lemma~\ref{lema1} and 
the recurrence formula (\ref{solution hypergeometric system}) that for every $p\ge1$ one has
$$
\left\|s_p\right\|\le \prod_{j=1}^{p}\left\|\left(\frac{m(j+\mu)}{m(j+\mu-1)}I-B\right)^{-1}\right\|\left\|A\right\|^{p}\left\|s_{0}\right\|\le \left\|s_{0}\right\|(\left\|A\right\|C)^p
$$
Normal convergence of the series provides holomorphy on some neighborhood of the origin of  $\sum_{p=0}^{\infty}s_pz^p$.
\end{proof}

\begin{remark} Observe from the proof of the previous result that condition (H2) only provides a sufficient condition for convergence. We have preferred to maintain such handeable hypotheses for the sake of simplicity. In the following result, another easier-to-handle sufficient condition is provided, not needing to compute an infinite number of inverse matrices, and can be directly checked from the sequence of moments.
\end{remark}

\vspace{0.3cm}

%\begin{corollary}
%Let $B\in\C^{n\times n}$ and $m=(m(p))_{p\ge0}$ satisfying (H1) and (H2). Then, $y(z)=z^{\mu}s_0$ is a solution of 
%$$z\partial_m y=By,$$
%where $s_0$ is an eigenvector associated to the eigenvalue $m(\mu)/m(\mu-1)$ of $B$.
%\end{corollary}

\begin{corollary}\label{coro1}
Under the hypotheses of Lemma~\ref{lema1}, assume $m(p+\mu)\neq0$ for all $p\ge0$ together with the sequence $(m(p+\mu-1)/m(p+\mu))_{p\ge1}$ being upper bounded. In addition to this, assume that 
\begin{equation}\label{e474}
\left\|B\right\|<\left|\frac{m(p+\mu)}{m(p+\mu-1)}\right|,
\end{equation}
for every $p\ge1$. Then, the system (\ref{eqhypsys}) admits a unique convergent solution in the form (\ref{e416}).
\end{corollary}
\begin{proof}
We observe that for every $p\ge1$ one has
$$\left\|\left(\frac{m(p+\mu)}{m(p+\mu-1)}I- B\right)^{-1}\right\|=\left|\frac{m(p+\mu-1)}{m(p+\mu)}\right|\left\|I-\left(\frac{m(p+\mu-1)}{m(p+\mu)}B\right)^{-1}\right\|.$$
In view of (\ref{e474}), the previous norm can be upper bounded by 
$$\left|\frac{m(p+\mu-1)}{m(p+\mu)}\right|\sum_{j\ge0}\left\|\frac{m(p+\mu-1)}{m(p+\mu)}B\right\|^{j},$$
which is a normally convergent series. The application of Corollary~\ref{coro0} concludes the result.
\end{proof}

We provide some examples of the previous result regarding the convergence of the formal solution of~(\ref{eqhypsys}). 

\begin{example}
Let $(m(p))_{p\ge0}$ be the sequence of Catalan numbers, i.e, $m(p)=\frac{(2p)!}{(p+1)!p!}$, for all $p\ge0$. We also fix $B\in \C^{2x2}$ given by
\begin{equation*}
    B=\begin{pmatrix}
        4-\frac{6}{\mu+1}&1\\
        0& 4-\frac{6}{\mu+1},
    \end{pmatrix}
\end{equation*}
where $\mu$ is a positive integer. Hypothesis (H1) holds, and therefore Lemma~\ref{lema1} guarantees the existence of a formal power series in the form (\ref{e416}) of 
$$z\partial_my=(zA+B)y,$$
for any choice of $A\in\C^{2\times 2}$. In addition to this, if one considers the matrix 1-norm, which is induced by the corresponding vector 1-norm, one has that
\begin{align*}
    \left\|\biggl(\bigl(4-\frac{6}{p+\mu+1}\bigr)I-B\biggr)^{-1}\right\| &= \frac{1}{36 p^2}\biggl| (\mu +1) (p+\mu +1)\biggr|  \biggl(\biggl| (\mu +1) (p+\mu +1)\biggr| +6 p\biggr)\\&\le \frac{1}{36} \biggl| (\mu +1) (\mu +2)\biggr|  \biggl(\biggl| (\mu +1) (\mu +2)\biggr| +6\biggr).
\end{align*}
Therefore, (H2) holds and Corollary~\ref{coro0} states the convergence of the formal solution. Indeed, the convergent solution is given by 
$$y(z) = \sum_{p=0}^\infty s_pz^{p+\mu},$$
with $s_0=(1,0)$, and 
$$s_{p}=\biggl(\bigl(4-\frac{6}{p+\mu+1}\bigr)I-B\biggr)^{-1}As_{p-1}, \quad p\geq 1.$$
\end{example}

\begin{example}
Let $B\in\C^{n\times n}$, for some integer $n\ge1$. Assume there exists $\mu\in\C$ with $\hbox{Re}(\mu)\ge 1$ and $C_1,C_2>0$ with
$$C_1\le \left\|B\right\|\le (m(p+\mu)/m(p+\mu-1))\le C_2,$$
for every $p\ge0$, and some submultiplicative matrix norm $\left\|\cdot\right\|$. Then, Corollary~\ref{coro0} states the convergence of the formal solution to 
$$z\partial_my=(zA+B)y,$$ 
in the form (\ref{e416}), for every $A\in\C^{n\times n}$.
\end{example}

% \color{red}
% \textbf{QUITAR}
% \begin{example}
% Let $m(z)$ be the sequence of Catalan numbers, i.e, the evaluation of $m(z) = \frac{4^{z}}{\sqrt{\pi}}\frac{\Gamma(z+\frac{1}{2})}{\Gamma(z+2)}$ at the natural numbers, and let $B\in \C^{2x2}$ be given by
% \begin{equation*}
%     B=\begin{pmatrix}
%         4-\frac{6}{\mu+1}&1\\
%         0& 4-\frac{6}{\mu+1},
%     \end{pmatrix}
% \end{equation*}
% where $\mu\in \C$ such that $Re(\mu)\ge 1$ so (H1) holds. Moreover taking the 1-norm, which is simply the maximum absolute column sum of the matrix, we have 
% \begin{align*}
%     \left\|\biggl(\bigl(4-\frac{6}{n+\mu+1}\bigr)I-B\biggr)^{-1}\right\| &= \frac{1}{36 n^2}\biggl| (\mu +1) (n+\mu +1)\biggr|  \biggl(\biggl| (\mu +1) (n+\mu +1)\biggr| +6 n\biggr)\\&\le \frac{1}{36} \biggl| (\mu +1) (\mu +2)\biggr|  \biggl(\biggl| (\mu +1) (\mu +2)\biggr| +6\biggr).
% \end{align*}
% In such a situation, the hypotheses (H1) and (H2) are satisfied. Therefore, according to Corollary \ref{coro0}, the formal series  $$y(z) = \sum_{n=0}^\infty s_nz^{n+\mu},$$
% with $$Bs_0=\bigl(4-\frac{6}{\mu+1}\bigr)s_0$$ and $$s_{n}=\biggl(\bigl(4-\frac{6}{n+\mu+1}\bigr)I-B\biggr)^{-1}As_{n-1}, \quad n\geq 1,$$
% is the unique convergent solution to (\ref{eqhypsys}).
% \end{example}

\color{black}

As a result Lemma~\ref{lema1}, Corollary~\ref{coro0} and Corollary~\ref{coro1}, the following result holds.

\begin{theorem}\label{teopral1}
Consider the linear system of moment differential equations (\ref{eqhypsys}).
\begin{itemize}
\item Let $1\le k\le n$. The system (\ref{eqhypsys}) admits $k$ linearly independent formal solutions of the form (\ref{e416}) if and only if $B$ admits $k$ linearly independent eigenvectors $s_j$, $1\le j\le k$, such that for every $1\le j\le k$ there exists $\mu_j\in\C$ with $\hbox{Re}(\mu)\ge1$ and $\frac{m(\mu_j)}{m(\mu_j-1)}\in\hbox{spec}(B)$, with $\frac{m(p+\mu_j)}{m(p+\mu_j-1)}\not\in\hbox{spec}(B)$ for every $p\ge1$. 
\item Let $1\le k\le n$. If the following assumptions hold:
\begin{itemize}
\item[(a)] $B$ admits $k$ linearly independent eigenvectors $s_j$, $1\le j\le k$, such that for every $1\le j\le k$ there exists $\mu_j\in\C$ with $\hbox{Re}(\mu_j)\ge1$ and $\frac{m(\mu_j)}{m(\mu_j-1)}\in\hbox{spec}(B)$, with $\frac{m(p+\mu_j)}{m(p+\mu_j-1)}\not\in\hbox{spec}(B)$ for every $p\ge1$.
\item[(b)] There exists $C>0$ such that for every $1\le j\le k$ and all $p\ge1$
$$\left\|\left(\frac{m(p+\mu_j)}{m(p+\mu_j-1)}I- B\right)^{-1}\right\|\le C,$$
\end{itemize}
then the system (\ref{eqhypsys}) admits $k$ linearly independent analytic solutions of the form (\ref{e416}).
\end{itemize}
\end{theorem}

\begin{example}
    Let $m(z) = \Gamma(1+\frac{1}{z})$ and let $B\in \C^{2x2}$ be given by 
    \begin{equation*}
     B= 
    \begin{pmatrix}
        \frac{\Gamma(1+\frac{1}{\mu_1})}{\Gamma(1+\frac{1}{\mu_1-1})}&0\\
        1&\frac{\Gamma(1+\frac{1}{\mu_2})}{\Gamma(1+\frac{1}{\mu_2-1})}
        
    \end{pmatrix},\quad \mu_1,\mu_2\in\C
      \end{equation*}
      where $\mu_1- \mu_2 \notin \Z$ and $Re(\mu_1),Re(\mu_2)\ge 1$. In such a situation, $B$ admits $2$ linearly independent eigenvectors
      \begin{equation*}
    s_{1,0} = \begin{pmatrix}
    \frac{\Gamma(1+\frac{1}{\mu_1})}{\Gamma(1+\frac{1}{\mu_1-1})}- \frac{\Gamma(1+\frac{1}{\mu_2})}{\Gamma(1+\frac{1}{\mu_2-1})} \\
    1
    \end{pmatrix},\quad s_{2,0}= \begin{pmatrix}
    0\\1
    \end{pmatrix}
\end{equation*}
Such that $Bs_{j,0} = \frac{\Gamma(1+\frac{1}{\mu_j})}{\Gamma(1+\frac{1}{\mu_j-1})}s_{j,0}$ with $j=1,2$. Then,  for all $n\geq 1$, the system (\ref{eqhypsys}) admits two linearly independent power series solutions given by 
\begin{equation*}
    y_j(z) = \sum_{p=0}^\infty s_{j,p}z^{p+\mu_j}, \quad j=1,2
\end{equation*}
with 
$$s_{j,p} = \biggl(\frac{\Gamma(1+\frac{1}{p+\mu_j})}{\Gamma(1+\frac{1}{p+\mu_j-1})}I-B\biggr)^{-1}As_{j,p-1}.$$ 

If $\left\|\cdot\right\|$ is a submultiplicatively norm and $B\in\C^{2\times 2}$ are such that
$$\left\|B\right\|\le\inf_{j=1,2}\inf_{p\ge0}\left|\frac{\Gamma(1+\frac{1}{p+\mu_j})}{\Gamma(1+\frac{1}{p+\mu_j-1})}\right|,$$
then Corollary~\ref{coro1} guarantees convergence of the previous solutions. Observe that, in particular, $\left\|B\right\|\le 1$ is a necessary condition (compute the limit when $n \to\infty$).
\end{example}

%\color{blue}

%Moreover,  $y_{1,2}(z)$ define two linearly independent analytic solutions provided that 
%\begin{align*}
%    \left\|\biggl(\frac{m(n+\mu_j)}{m(n+\mu_j-1)}I-B\biggr)^{-1}\right\|&= \frac{1+\biggl|\frac{\Gamma(1+\frac{1}{\mu_2})}{\Gamma(1+\frac{1}{\mu_2-1})}-\frac{\Gamma(1+\frac{1}{n+\mu_j})}{\Gamma(1+\frac{1}{n+\mu_j-1})}\biggr|}{\biggl|\frac{\Gamma(1+\frac{1}{\mu_1})}{\Gamma(1+\frac{1}{\mu_1-1})}-\frac{\Gamma(1+\frac{1}{n+\mu_j})}{\Gamma(1+\frac{1}{n+\mu_j-1})}\biggr|\biggl|\frac{\Gamma(1+\frac{1}{\mu_2})}{\Gamma(1+\frac{1}{\mu_2-1})}-\frac{\Gamma(1+\frac{1}{n+\mu_j})}{\Gamma(1+\frac{1}{n+\mu_j-1})}\biggr|}\\ &\le\frac{1+\biggl|\frac{\Gamma(1+\frac{1}{\mu_2})}{\Gamma(1+\frac{1}{\mu_2-1})}-\frac{\Gamma(1+\frac{1}{1+\mu_j})}{\Gamma(1+\frac{1}{\mu_j})}\biggr|}{\biggl|\frac{\Gamma(1+\frac{1}{\mu_1})}{\Gamma(1+\frac{1}{\mu_1-1})}-\frac{\Gamma(1+\frac{1}{1+\mu_j})}{\Gamma(1+\frac{1}{\mu_j})}\biggr|\biggl|\frac{\Gamma(1+\frac{1}{\mu_2})}{\Gamma(1+\frac{1}{\mu_2-1})}-\frac{\Gamma(1+\frac{1}{1+\mu_j})}{\Gamma(1+\frac{1}{\mu_j})}\biggr|}. 
%\end{align*}
 
%\end{example}
%\color{black}

\begin{remark}
Observe that the classical case, i.e. $m=(p!)_{p\ge0}$ the moment derivative is reduced to the usual derivation. The previous result states the existence of a solution of the system when two eigenvalues of $B$ do not differ by a positive integer, i.e. $B$ and $B+pI$ share no eigenvalues for every $p\ge1$. In the literature, such systems are said to have good spectrum, or the matrix $B$ being non-resonant (see~\cite{balser}).
\end{remark}

\subsection{Functions of a matrix}\label{secfunmat}

In~\cite{lastra}, the first author describes the general solution to (\ref{eqhypsys}) in the case that $B\equiv 0$. In~\cite{lastraprisuelos}, the solution is written in the form of a generalized matrix exponential.

It the present situation, i.e. $B\neq 0$, the classical theory states that the matrix $z^{B}$ turns out to be a fundamental solution of the system (\ref{eqhypsys}) whenever $A\equiv 0$, and $k=n$, in Theorem~\ref{teopral1}. In a more general framework, one can also extend the definition of $z^{B}$ as follows:

\begin{definition}\label{defin1}
In the situation of Theorem~\ref{teopral1}, let us assume that $A\equiv 0$ in (\ref{eqhypsys}), and $k=n$, i.e. there exist $n$ linearly independent Floquet solutions. We define the matrix
$$z_m^{B}:=(C_1(z),\ldots,C_n(z)),$$
where $C_j(z)\in\C[[z]]^{n}$ is a column vector given by $s_{0,j}z^{\mu_j}$, where $s_{0,j}$ stands for the eigenvector of $B$ associated with the eigenvalue $m(\mu_j)/m(\mu_j-1)$ for every $1\le j\le n$. 
\end{definition}

\begin{example}
Let $B=\hbox{diag}(b_1,\ldots,b_n)\in\C^{n\times n}$. The matrix $z_m^{B}$, if it exists, is given by
$$z_m^{B}=\hbox{diag}(z^{\mu_1},\ldots,z^{\mu_n}),$$
where $\mu_j\in\C$ is such that $m(\mu_j)/m(\mu_j-1)=b_j$, $1\le j\le n$. The existence of such matrix is conditioned on the existence of such $\mu_j$, $1\le j\le n$, which is always guaranteed whenever $m=(p!)_{p\ge0}$, i.e. $\partial_m$ is the classical derivative, with $m(\mu_j)/m(\mu_j-1)=\Gamma(\mu_j)/\Gamma(\mu_j-1)=\mu_j=b_j$. 
\end{example}

The previous definition of $z_m^{B}$ constitutes a first approach to monodromy results for systems of moment differential equations. For a background on the classical theory, we refer to~\cite{haraoka}. As a matter of fact, Theorem 4.5. (i) in~\cite{haraoka} describes the fundamental matrix solution associated to a linear system of differential equations of the form $y'=B(z)y$, with $B=B(z)$ being an $n\times n$ matrix with meromorphic coefficients.

\subsection{A refinement of the results}\label{sec35}

Notice that the assumption $\hbox{Re}(\mu)>1$ in (H1) can be obviated. The construction of a formal (resp. a convergent) solution follows analogous arguments under the following assumptions. Let us define the following assumptions:

\begin{itemize}
    \item [(H1)'] There exist $\mu \in \C$ and a nonnegative integer $N$ such that $Re(\mu)+N \geq 1$ together with $\frac{m(\mu+N)}{m(\mu+N-1)}\in spec(B)$ and $\frac{m(p+\mu+N)}{m(p+\mu+N-1)}\notin spec(B)$, for any positive integer $p$. 
    \item [(H2)'] There exists $C'>0$ such that
    $$\left\|\left(\frac{m(p+\mu+N)}{m(p+\mu+N-1)}I- B\right)^{-1}\right\|\le C',$$
for every $p\ge1$.
\end{itemize}

A generalized version of Lemma~\ref{lema1} and Corollary~\ref{coro0} read as follows. Their proof can follow analogous arguments as that of their related previous results, so we omit them.

\begin{lemma}
Let (\ref{eqhypsys}) be a linear system of moment differential equations such that hypothesis (H1)' holds. Assume that $s_{0,N}$ is an eigenvector associated to the eigenvalue $\frac{m(\mu+N)}{m(\mu+N-1)}$ of $B$. Then, the formal power series 
\begin{equation}\label{e419}
y(z)=\sum_{p=0}^{\infty}s_{p,N}z^{p+\mu+N},
\end{equation}
where
\begin{align}\label{solution hypergeometric system b}
s_{p,N}=\left(\frac{m(p+\mu+N)}{m(p+\mu+N-1)}I-B\right)^{-1}As_{p-1,N},\quad p\geq 1.
\end{align}
is a formal solution to (\ref{eqhypsys}).
\end{lemma}

\begin{lemma}
    Let the assumptions (H1)' and (H2)' hold.  Let $s_{0,N}$ be an eigenvector of the eigenvalue $\frac{m(\mu+N)}{m(\mu+N-1)}$ of $B$. Then, the power series (\ref{e419}) defined by (\ref{solution hypergeometric system b}) is the unique convergent solution to (\ref{eqhypsys}) of the form (\ref{e419}).
\end{lemma}
\begin{proof}
This proof follows a similar approach to the proof of \ref{lema1} by considering the formal power series \eqref{e419} as a solution to (\ref{eqhypsys}).
% Plugging the formal series \eqref{e419} into (\ref{eqhypsys}), we get the recurrence 
% \begin{equation*}
%    Bs_{0,N} =  \frac{m(\mu+N)}{m(\mu+N-1)}s_{0,N}, \quad n = 0,
% \end{equation*}
% and
% \begin{equation*}
%     \frac{m(\mu+N+n)}{m(\mu+N+n-1)}s_{n,N} = As_{n-1,N} + Bs_{n,N}, \quad \forall n\geq 1.
% \end{equation*}
% Indeed
% \begin{align*}
%     z\partial_m y(z;N) &= z\sum_{n=0}^\infty s_{n,N}\frac{m(\mu+N+n)}{m(\mu+N+n-1)}z^{n+N+\mu-1} = \sum_{n=0}^\infty s_{n,N}\frac{m(\mu+N+n)}{m(\mu+N+n-1)}z^{n+N+\mu} \\ & = Bs_{0,N}z^{\mu +N}+ \sum_{n=1}^\infty s_{n,N}\frac{m(\mu+N+n)}{m(\mu+N+n-1)}z^{n+N+\mu} =  A\sum_{n=1}^\infty s_{n-1,N}z^{n+\mu+N}+B\sum_{n=0}^\infty s_{n,N}z^{n+\mu+N} \\ & = A\sum_{n=0}^\infty s_{n,N}z^{n+\mu+N+1}+B\sum_{n=0}^\infty s_{n,N}z^{n+\mu+N} = \bigl(zA+B\bigr)\sum_{n=0}^\infty s_{n,N}z^{n+\mu+N} = \bigl(zA+B\bigr)y(z;N)
% \end{align*}
% Moreover, from the general recurrence and assumption (H2)'
% \begin{equation*}
%     \left\|s_{n,N}\right\| \le \prod_{j=1}^n \left\|\biggl(\frac{m(j+\mu+N)}{m(j+\mu+N-1)}\biggr)^{-1}\right\|\left\|A\right\|^n\left\|s_{0,N}\right\|\le \left\|s_{0,N}\right\|(\left\|A\right\|C')^n
% \end{equation*}
\end{proof}

\subsection{Some important applications of the theory}\label{secapap}

We particularize the previous results to two specific moment sequences due to their significance in applications. The first result describes the $q$-Gevrey setting. Given $q>1$, $q$-Gevrey sequence  $([p]_q^{!})_{p\ge0}$ is not a strongly regular sequence. However, one can construct a pair of kernel-like functions together with their associated sequence of moments, as described in Section~\ref{sec2}. Therefore, the results obtained remain valid in this context.    %the $q$-Gevrey sequence, $m_q=([p]_q!)_{p\ge 0}$, defined in Section 3, which is no longer a strongly regular sequence since (mg) condition is not satisfied.

First, we consider systems of $q-$difference equations. We recall that given $q>1$, the particularization of a moment derivative to the $q-$Gevrey sequence $m_q=([p]_q^{!})_{p\ge0}$ determines Jackson $q-$derivative $D_q$:
$$D_q(f)(z)=\frac{f(qz)-f(z)}{qz-z}.$$
On the other hand, $q-$Gamma function is defined by
$$\Gamma_{q}(z)=\frac{(q^{-1};q^{-1})_{\infty}}{(q^{-z};q^{-1})_{\infty}}(q-1)^{1-z}q^{\frac{z(z-1)}{2}},$$
for every $z\in\{z\in\C:\hbox{Re}(z)>0\}$. For all $\omega\in\C\setminus(-\infty,0]$, we write $z^\omega$ for $e^{\omega\log(z)}$, where $\log(\cdot)$ is the principal branch of the logarithm, and
$$ (\alpha;q^{-1})_{\infty}=\prod_{p\ge0}(1-\frac{\alpha}{q^{p}}),\quad \alpha\in\C.$$
It turns out that $\Gamma_q$ is holomorphic on $\{z\in\C:\hbox{Re}(z)>0\}$ with 
$$\Gamma_q(p+1)=[p]_q^{!},\quad p\ge0$$
(see~\cite{as,zhang}, for example). Therefore, the previous results can be applied when considering the moment function to be $\Gamma_q(\cdot)$.

\begin{corollary}
    Let $q\in (1,\infty)$ and consider $m_q=([p]_q!)_{p\ge 0}$, and $A,B\in\C^{n\times n}$. Assume that the hypothesis (H1) holds. In particular, there exists $\mu\in\C$ with $\hbox{Re}(\mu)>1$ such that 
$$\frac{\Gamma_q(\mu)}{\Gamma_q(\mu-1)}=\frac{q^{\mu-1}-1}{q-1}=[\mu]_q$$ 
is an eigenvalue of $B$.
    
   % Let $\mu\in\C$ such that $Re(\mu)\ge 1$ and $[\mu]_q\in spec(B)$ while $[\mu+n]_q\notin spec(B)$, for every $n\in N$. 
    
Let $s_{0,q}$ be an eigenvector associated to the eigenvalue $[\mu]_q$. Then, the formal power series 
   \begin{equation*}
       y(z;q) = \sum_{p=0}^\infty s_{p,q}z^{p+\mu},
   \end{equation*}
   with 
   $$s_{p,q} = \biggl([\mu+p]_qI-B\biggr)^{-1}As_{p-1,q},$$
   is the unique convergent solution of the $q-$hypergeometric system 
   \begin{equation*}
       zD_qy(z;q) = (zA+B)y(z;q).
   \end{equation*}
	in the form of a formal Floquet series. In addition to this, if $\left\|\cdot\right\|$ is a submultiplicative matrix norm, and (H2) holds for some $C>0$, i.e.
	$$\sup_{p\ge 1}\left\|[p+\mu]_qI-B)^{-1}\right\|\le C,$$
	then, the previous Floquet formal solution is convergent.
\end{corollary}

As a second application, we describe the main result of the section in the framework of Caputo fractional derivatives, quite related to moment derivation when considering the sequence of moments $m_{\alpha}$, for some positive $\alpha$.

\begin{corollary}
    Let $\alpha>0$ and consider $m_\alpha=\bigl(\Gamma(1+\frac{p}{\alpha})\bigr)_{p\ge 0}$. Let $A,B\in\C^{n\times n}$, and assume that the hypothesis (H1) holds, for some $\mu\in\C$ with $\hbox{Re}(\mu)\ge1$. Let $s_{0,\alpha}$ is an eigenvector associated to the eigenvalue $\frac{\Gamma(1+\frac{\mu}{\alpha})}{\Gamma(1+\frac{\mu-1}{\alpha})}$ of $B$. Then, the formal power series 
   $$
       y(z;\alpha) = \sum_{p=0}^\infty s_{p,q}z^{(p+\mu)/\alpha},
   $$
   with 
   $$s_{p,q} = \biggl(\frac{\Gamma(1+\frac{\mu+p}{\alpha})}{\Gamma(1+\frac{\mu+p-1}{\alpha})}I-B\biggr)^{-1}As_{p-1,q},$$
   is the unique formal solution of the fractional hypergeometric system 
   \begin{equation*}
       z^{1/\alpha} \prescript{C}{}{D_z}^{1/\alpha}y(z) = (z^{1/\alpha}A+B)y(z)
   \end{equation*}
	in the form of a Floquet series. In addition to this, if if $\left\|\cdot\right\|$ is a submultiplicative matrix norm, and (H2) holds for some $C>0$, i.e.
	$$\sup_{p\ge 1}\left\|(\frac{\Gamma(1/\alpha)}{B(1+\frac{\mu-1}{\alpha},\frac{1}{\alpha})}I-B)^{-1}\right\|\le C,$$
	then, the previous Floquet formal solution is convergent.
\end{corollary}

\section{Confluent hypergeometric systems: structure of the solution}

This section is devoted to prove several results. A first technical lemma (Lemma~\ref{lematec}) will be applied in the rest of the study. From the point of view of the authors, it deserves particular attention for its importance, which relies on the first application of a change of variable in the framework of moment differential equations. In general, moment differentiation does not satisfy a Leibniz-like rule, which makes it hard to find such techniques available in the general framework.

A second result, Proposition~\ref{prop1}, describes the structure of the solutions, in the spirit of the classical theory of systems of differential equations. We conclude the section by showing the explicit form of solutions of any system in dimension $n=2$. The previous results allow to reduce such situations to those considered in Section~\ref{secdiag} and Section~\ref{secnodiag}.

In the whole section, we assume $m=(m(p))_{p\ge0}$ is a sequence of positive real numbers determined by the restriction of some function defined on $\{z\in\C:\hbox{Re}(z)\ge0\}$. In particular, one can think about the sequence of moments associated to a pair of kernel functions $(e,E)$ for generalized summability, or the $q$-Gevrey sequence.

\begin{lemma}\label{lematec}
Let us consider $A,B\in\C^{n\times n}$, for some $n\ge1$. Let $\lambda\in\hbox{spec}(A)$. We  assume the matrices $A$ and $B$ commute. We also take for granted Assumption $(H1)$. Let $\lambda\in\C$. Then, $y=y(z)$ is a solution of the linear system of moment differential equations
\begin{equation}\label{e496}
z\partial_m y=(zA+B)y
\end{equation}
in the form (\ref{e416}) if and only if
\begin{align}\label{eq:change of variables}
y(z)=h(z)\widetilde{y}(z),
\end{align}
where $\tilde{y}$ is a solution of
\begin{align}\label{eq: hypergeometric det=0}
z\partial_m\widetilde{y}(z)=(z(A-\lambda I)+B)\widetilde{y}(z).
\end{align}
Here, $h(z)=\sum_{p=0}^{\infty}h_pz^p\in\C^{n\times n}[[z]]$, is defined by $h_0=I$ and
\begin{multline*}
\left(\frac{m(p+\mu)}{m(p+\mu-1)}-\frac{m(\mu)}{m(\mu-1)}\right) h_{p}\\
=\sum_{j=0}^{p-1}h_j\left(\lambda \widehat{s}_{p-j-1}-\left(\frac{m(p+\mu)}{m(p+\mu-1)}-\frac{m(p+\mu-j)}{m(p+\mu-j-1)}\right)\widehat{s}_{p-j}\right),\,\textrm{ for all}\quad p\geq 1;
\end{multline*}
where $\widehat{s}_{0}=I$ and
$$ \widehat{s}_{k}:=\prod_{i=0}^{k-1}\left(\left(\frac{m(k+\mu-i)}{m(k+\mu-i-1)}I-B\right)^{-1}(A-\lambda I)\right), \,k\geq 1.$$
\end{lemma}

\begin{proof}
We prove the  change of variable  \eqref{eq:change of variables} transforms the system  (\ref{e496})  into \eqref{eq: hypergeometric det=0}.

Some first remarks are at hand. Hypothesis $(H1)$ guarantees the existence of inverse for the matrix
$$\frac{m(k+\mu-i)}{m(k+\mu-i-1)}I-B$$
for every $0\le i\le k-1$ and all $k\ge1$. Also, $(H1)$ determines that $\frac{m(p+\mu)}{m(p+\mu-1)}\neq\frac{m(\mu)}{m(\mu-1)}$ for every $p\ge1$. Therefore, the recurrence determining the sequence $(h_p)_{p\ge0}$ is well-defined.

Let $y(z)$ in the form \eqref{eq:change of variables}. The result follows from the following fact:
\begin{align}\label{eq:lemma chage of variables}
\partial_my(z)=h(z)\partial_m\widetilde{y}(z)+\lambda h(z)\widetilde{y}(z).
\end{align}

\begin{lemma}\label{lema530}
$h$ commutes with $A$ and $B$.
\end{lemma}
\begin{proof}
From the definition of $h$, one has that $h$ commutes with $A$ and $B$ if and only if $\widehat{s}_{j}$ commutes with $A$ and $B$ for every $j\ge0$. From te definition of $\widehat{s}_{j}$, the statement is then reduced to check that $(\frac{m(\ell+\mu)}{m(\ell+\mu-1)}I-B)^{-1}$ and $A$ commute with $A$ and $B$, which is a consequence of $A$ and $B$ being matrices which commute.
\end{proof}

First, assume that (\ref{eq:lemma chage of variables}) holds. Then, we observe that (\ref{e496}) implies
$$(zA+B)y=zh(z)\partial_m\widetilde{y}+\lambda z h(z)\widetilde{y},$$
i.e.
$$z\partial_m\widetilde{y}=h(z)^{-1}(zA+B)h(z)h(z)^{-1}y-\lambda z\widetilde{y}.$$
Recall that $h$ is invertible due to $h_0=I$. Lemma~\ref{lema530} allows to conclude that $\widetilde{y}$ is a solution of (\ref{eq: hypergeometric det=0}). The converse statement follows in the same way.

It only rests to prove the equality (\ref{eq:lemma chage of variables}). 

Let us write $\widetilde{y}(z)=\sum_{p=0}^{\infty}\widehat{s}_{p}\widetilde{s}_0 z^{p+\mu}$, for unknown $\widehat{s}_{p}$ for $p\ge0$ and where $\widetilde{s}_0=s_0$ is as above. Let us write $\widetilde{s}_p:=\widehat{s}_{p}\widetilde{s}_0$ for every $p\ge0$. In order to compute $\partial_m y(z)$, note that the application of Cauchy product in the relation (\ref{eq:change of variables}) leads us to
\begin{align*}
y(z)=\sum_{p=0}^{\infty}z^{p+\mu}\sum_{j=0}^ph_j\widetilde{s}_{p-j}.
\end{align*}
Therefore,
$$
\partial_m y(z)=\sum_{p=0}^{\infty}\frac{m(p+\mu)}{m(p+\mu-1)}z^{p+\mu-1}\sum_{j=0}^ph_j\widetilde{s}_{p-j}.
$$
Moreover, one has that
\begin{align*}
h(z)\partial_m\widetilde{y}(z)&=\left(\sum_{p=0}^{\infty}h_pz^p\right)\sum_{p=0}^{\infty}\widetilde{s}_p\frac{m(p+\mu)}{m(p+\mu-1)}z^{p+\mu-1}\\
&=\sum_{p=0}^{\infty}z^{p+\mu-1}\sum_{j=0}^{p}h_j\widetilde{s}_{p-j}\frac{m(p-j+\mu)}{m(p-j+\mu-1)}.
\end{align*}
Consequently, for all $\mu \geq 0$,
\begin{align*}
\partial_m y(z)-h(z)\partial_m\widetilde{y}(z)&=
\sum_{p=1}^{\infty}z^{p+\mu-1}\sum_{j=0}^p\left(\frac{m(p+\mu)}{m(p+\mu-1)}-\frac{m(p-j+\mu)}{m(p-j+\mu-1)}\right)h_j\widetilde{s}_{p-j}
\\
&
=
\sum_{p=0}^{\infty}z^{p+\mu}\sum_{j=0}^{p+1}\left(\frac{m(p+1+\mu)}{m(p+\mu)}-\frac{m(p+1-j+\mu)}{m(p-j+\mu)}\right)h_j\widetilde{s}_{p+1-j}
\\
&
=
\sum_{p=0}^{\infty}z^{p+\mu}\sum_{j=1}^{p+1}\left(\frac{m(p+1+\mu)}{m(p+\mu)}-\frac{m(p+1-j+\mu)}{m(p-j+\mu)}\right)h_j\widetilde{s}_{p+1-j}
\\
&
=
\sum_{p=0}^{\infty}z^{p+\mu}\sum_{j=0}^{p}\left(\frac{m(p+1+\mu)}{m(p+\mu)}-\frac{m(p-j+\mu)}{m(p-j+\mu-1)}\right)h_{j+1}\widetilde{s}_{p-j}.
\end{align*}
On the other hand, one can compute $\lambda h(z)\widetilde{y}(z)$ analogously:
\begin{align*}
\lambda h(z)\widetilde{y}(z)=\sum_{p=0}^{\infty}z^{p+\mu}\sum_{j=0}^p\lambda h_j\widetilde{s}_{p-j}.
\end{align*}
Notice that
\begin{multline*}
\sum_{j=0}^{p}\left(\left(\frac{m(p+1+\mu)}{m(p+\mu)}-\frac{m(p-j+\mu)}{m(p-j+\mu-1)}\right)h_{j+1}-
\lambda h_j\right)\widetilde{s}_{p-j}\\
=
\sum_{j=0}^{p}\left(\frac{m(p+1+\mu)}{m(p+\mu)}-\frac{m(p-j+\mu)}{m(p-j+\mu-1)}\right)h_{j+1}\widehat{s}_{p-j}\widetilde{s}_0-\sum_{j=0}^{p}
\lambda h_j\widehat{s}_{p-j}\widetilde{s}_0\hfill\\
=
\sum_{j=0}^{p}\left(\frac{m(p+1+\mu)}{m(p+\mu)}-\frac{m(p-j+\mu)}{m(p-j+\mu-1)}\right)h_{j+1}\widehat{s}_{p-j}\widetilde{s}_0-\sum_{j=1}^{p}
\lambda h_j\widehat{s}_{p-j}\widetilde{s}_0-\lambda h_0\widehat{s}_{p}\widetilde{s}_0
\hfill\\
=
\sum_{j=0}^{p}\left(\frac{m(p+1+\mu)}{m(p+\mu)}-\frac{m(p-j+\mu)}{m(p-j+\mu-1)}\right)h_{j+1}\widehat{s}_{p-j}\widetilde{s}_0-\sum_{j=0}^{p-1}
\lambda h_{j+1}\widehat{s}_{p-1-j}\widetilde{s}_0-\lambda h_0\widehat{s}_{p}\widetilde{s}_0
\hfill\\
=\left(\frac{m(p+1+\mu)}{m(p+\mu)}-\frac{m(\mu)}{m(\mu-1)}\right)h_{p+1}\widehat{s}_{0}\widetilde{s}_0\hfill \\
\hfill+
\sum_{j=0}^{p-1}h_{j+1}\left(\left(\frac{m(p+1+\mu)}{m(p+\mu)}-\frac{m(p-j+\mu)}{m(p-j+\mu-1)}\right)\widehat{s}_{p-j}-
\lambda \widehat{s}_{p-1-j}\right)\widetilde{s}_0-\lambda h_0\widehat{s}_{p}\widetilde{s}_0
\\
=\left(\frac{m(p+1+\mu)}{m(p+\mu)}-\frac{m(\mu)}{m(\mu-1)}\right)h_{p+1}\widetilde{s}_0\hfill
\\
\hfill+
\sum_{j=1}^{p}h_{j}\left(\left(\frac{m(p+1+\mu)}{m(p+\mu)}-\frac{m(p+1-j+\mu)}{m(p-j+\mu)}\right)\widehat{s}_{p+1-j}-
\lambda \widehat{s}_{p-j}\right)\widetilde{s}_0-\lambda h_0\widehat{s}_{p}\widetilde{s}_0
\\
=\left(\frac{m(p+1+\mu)}{m(p+\mu)}-\frac{m(\mu)}{m(\mu-1)}\right)h_{p+1}\widetilde{s}_0
+
\sum_{j=0}^{p}h_{j}\left(\left(\frac{m(p+1+\mu)}{m(p+\mu)}-\frac{m(p+1-j+\mu)}{m(p-j+\mu)}\right)\widehat{s}_{p+1-j}-
\lambda \widehat{s}_{p-j}\right)\widetilde{s}_0.
\end{multline*}
This entails that
$$\partial_m y(z)-h(z)\partial_m\widetilde{y}(z)-\lambda h(z)\widetilde{y}(z)$$
equals
\begin{multline*}
\sum_{p=0}^{\infty}z^{p+\mu}\sum_{j=0}^{p}\left(\left(\frac{m(p+1+\mu)}{m(p+\mu)}-\frac{m(p-j+\mu)}{m(p-j+\mu-1)}\right)h_{j+1}-
\lambda h_j\right)\widetilde{s}_{p-j}\\
=\sum_{p=0}^{\infty}z^{p+\mu}\left[\sum_{j=0}^p h_j\left(\left(\frac{m(p+\mu+1)}{m(p+\mu)}-\frac{m(p+\mu+1-j)}{m(p+\mu-j)}\right)\widehat{s}_{p+1-j}-\lambda\widehat{s}_{p-j}\right)\right.\\
\left.\qquad\qquad\qquad+\left(\frac{m(p+\mu+1)}{m(p+\mu)}-\frac{m(\mu)}{m(\mu-1)}\right)h_{p+1}\right]\widetilde{s}_0\\
=\sum_{p=0}^{\infty}z^{p+\mu}\left[-\left(\frac{m(p+\mu+1)}{m(p+\mu)}-\frac{m(\mu)}{m(\mu-1)}\right)h_{p+1}+\left(\frac{m(p+\mu+1)}{m(p+\mu)}-\frac{m(\mu)}{m(\mu-1)}\right)h_{p+1}\right]\widetilde{s}_0=0,
\end{multline*}
which completes the proof.
\end{proof}

\begin{remark} Observe that, under the hypotheses of the previous result, $\widetilde{y}(z)=\sum_{p=0}^{\infty}\widehat{s}_{p}\widetilde{s}_0 z^{p+\mu}$, where $\tilde{s}_0=s_0$ is an eigenvector associated to the eigenvalue $m(\mu)/m(\mu-1)$ of $B$.
\end{remark}

\begin{remark} Observe that the previous lemma deals with formal solutions of systems of moment differential equations, and $h\in\C^{n\times n}[[z]]$ might also have zero radius of convergence.
\end{remark}

\begin{corollary}
Assume the hypotheses of Lemma~\ref{lematec} hold, together with Assumption $(H2)$. Then, $y$, $\widetilde{y}$, define convergent solutions of the systems (\ref{e496}) and (\ref{eq: hypergeometric det=0}), respectively.
\end{corollary}
\begin{proof}
The series $y=y(z)$ is convergent from Corollary~\ref{coro1}. Observe that 
$$\widetilde{y}(z)=\sum_{p=0}^{\infty}\widehat{s}_{p}\widetilde{s}_0 z^{p+\mu}$$
with
$$\left\|\widehat{s}_p\right\|\le \prod_{i=0}^{p-1}\left(\left\|\left(\frac{m(p+\mu-i)}{m(p+\mu-i-1)}I-B\right)^{-1}\right\|\left\|A-\lambda I\right\|\right)\le (C\left\|A-\lambda I\right\|)^{p},$$
for every $p\ge1$. The result follows from here.
\end{proof}

%%%
\begin{remark}
Note that in the classical setting, that is, when $m(p)=p!$ for $p\ge0$, we have that
$$
h_p=\frac{1}{p}\left(\sum_{j=0}^{p-1}\left(\lambda\widehat{s}_{p-j-1}-j\widehat{s}_{p-j}\right)h_j\right).
$$
Then, by induction it is straight to check that
$
h_p=\frac{\lambda^p}{p!}I. 
$ for every $p\ge0$. Hence, the change of variables in \eqref{eq:change of variables} leads us to
$$
y(z)=e^{\lambda z}\widetilde{y}(z),
$$
as expected.
\end{remark}

The following second simplification can also be achieved.

\begin{proposition}\label{prop1}
Let $m$ be a sequence of positive real numbers, and $A,B\in\C^{n\times n}$. Then, the $n\times n$ matrix $Y$ is a solution of (\ref{eqhypsys}) if and only if $P^{-1}Y$ is a solution of the system
$$z\partial_my=(z\tilde{A}+J)y,$$
where $B=PJP^{-1}$ is Jordan decomposition of the matrix $B$, and $\tilde{A}\in\C^{n\times n}$ with $\hbox{det}(\tilde{A})=\hbox{det}(A)$. 
\end{proposition}
\begin{proof}
Let $B=PJP^{-1}$ be Jordan decomposition of $B$, for some regular matrix $P$. Then, one has that 
$$z\partial_mY=(zAP+PJ)P^{-1}Y,$$
i.e.
$$z\partial_mY=P(zP^{-1}AP+J)P^{-1}Y.$$
Therefore,
$$P^{-1}z\partial_mY=(zP^{-1}AP+J)P^{-1}Y.$$
Linearity properties of $\partial_m$ lead us to
$$z\partial_m(P^{-1}Y)=(zP^{-1}AP+J)P^{-1}Y,$$
which allows us to conclude the result with $\tilde{A}=P^{-1}AP$. Observe that $\hbox{det}(A)=\hbox{det}(\tilde{A})$.
\end{proof}

\begin{example}
In the situation of Theorem~\ref{teopral1}, let us assume that $A\equiv0$ in (\ref{eqhypsys}), and $k=n$, i.e. there exist $n$ linearly independent Floquet solutions, $z_m^{B}$ (see Definition~\ref{defin1}). From Proposition~\ref{prop1} one gets that if $B=PJP^{-1}$ is the Jordan canonical decomposition of $B$, then the matrix
$$Y(z)=Pz_m^{J}$$
is a solution of (\ref{eqhypsys}).
\end{example}

The previous example is the moment version of Corollary 4.6,~\cite{haraoka}, in the classical framework, which is proved there under more generality on the matrix defining the main problem under study.

In the next sections we shall compute Floquet solution of (\ref{eqhypsys}) in the case of systems of dimention 2 in different situations. In view of Proposition~\ref{prop1}, $B$ can be chosen as a diagonal matrix or a matrix in Jordan form. Additionally, Lemma~\ref{lematec} allows us to assume that $\textrm{det}(A)=0$, by considering $\lambda\in\hbox{spec}(A)$ in Lemma~\ref{lematec} in some particular situations.

\subsection{Planar systems. $B$ diagonal}\label{secdiag}

We consider  the equation (\ref{eqhypsys}) with
$$
A:=\left(
\begin{matrix}
a&b
\\
c&d
\end{matrix}
 \right)\quad \textrm{and}	\quad  B:=\left(
\begin{matrix}
\frac{m(\mu_1)}{m(\mu_{1}-1)}&0
\\
0&\frac{m(\mu_2)}{m(\mu_2-1)}
\end{matrix}
 \right),
$$ where $Re(\mu_1), Re(\mu_2) \geq 1$ are such that 
$$\frac{m(\mu_k)}{m(\mu_{k}-1)}\notin\left\{\frac{m(p+\mu_j)}{m(p+\mu_j-1)}:j\in\{1,2\},p\ge1\right\},\hbox{ for }k=1,2.$$
We observe that hypothesis (H1) holds for $\mu=\mu_1$ and also for $\mu=\mu_2$.

It is possible to apply Lemma~\ref{lema1} twice to attain two formal solutions in the form of a Floquet power series. Let us denote 
$$
s_p:=\left(\begin{matrix}
f_p
\\
g_p
\end{matrix}   \right), 
$$
for all $p\ge0$. We observe that $f_0=1$, $g_0=0$ determines an eigenvector associated with the eigenvalue $\frac{m(\mu_1)}{m(\mu_{1}-1)}$. The recursion for $s_p$ in \eqref{solution hypergeometric system} gives us, for all $p\geq 1$ and $\mu = \mu_1$,
\begin{align}\label{particular solution diagonal case}
\left(\frac{m(p+\mu_1)}{m(p+\mu_1-1)}-\frac{m(\mu_1)}{m(\mu_{1}-1)}\right)f_p=af_{p-1}+bg_{p-1},
\\
\left(\frac{m(p+\mu_1)}{m(p+\mu_1-1)}-\frac{m(\mu_2)}{m(\mu_{2}-1)}\right)g_p=cf_{p-1}+dg_{p-1}.\nonumber
\end{align}
Therefore, for $\beta:=\frac{m(\mu_1)}{m(\mu_{1}-1)}-\frac{m(\mu_2)}{m(\mu_{2}-1)}$ and $\lambda=a+d$, we have for all $p\geq 1$,
\begin{align*}
&\left(\frac{m(p+1+\mu_1)}{m(p+\mu_1)}-\frac{m(\mu_1)}{m(\mu_{1}-1)}\right)\left(\frac{m(p+\mu_1)}{m(p+\mu_1-1)}-\frac{m(\mu_2)}{m(\mu_{2}-1)}\right)f_{p+1}
\\=&
a\left(\frac{m(p+\mu_1)}{m(p+\mu_1-1)}-\frac{m(\mu_2)}{m(\mu_{2}-1)}\right)f_{p}+b\left(\frac{m(p+\mu_1)}{m(p+\mu_1-1)}-\frac{m(\mu_2)}{m(\mu_{2}-1)}\right)g_{p}
\\&=
a\left(\frac{m(p+\mu_1)}{m(p+\mu_1-1)}-\frac{m(\mu_2)}{m(\mu_{2}-1)}\right)f_{p}+b(cf_{p-1}+dg_{p-1})
\\&=
a\left(\frac{m(p+\mu_1)}{m(p+\mu_1-1)}-\frac{m(\mu_2)}{m(\mu_{2}-1)}\right)f_{p}+bcf_{p-1}+d\left(\left(\frac{m(p+\mu_1)}{m(p+\mu_1-1)}-\frac{m(\mu_1)}{m(\mu_{1}-1)}\right)f_p-af_{p-1}\right)
\\&=
\left(\lambda\left(\frac{m(p+\mu_1)}{m(p+\mu_1-1)}-\frac{m(\mu_1)}{m(\mu_1-1)}\right)+a\beta\right)f_{p}+(bc-da)f_{p-1}.
\end{align*}

We have
\begin{align*}
g_0=0,\quad g_1=\frac{c}{\frac{m(\mu_1+1)}{m(\mu_1)}-\frac{m(\mu_1)}{m(\mu_1-1)}+\beta},\quad f_0=1\quad \textrm{and}\quad  f_1=\frac{a}{\frac{m(\mu_1+1)}{m(\mu_1)}-\frac{m(\mu_1)}{m(\mu_1-1)}}.
\end{align*}

If $A$ and $B$ do not commute, then the previous recursion determines all the values of $f_p,g_p$, for $p\ge0$. Otherwise, the computations become more handly, and in virtue of Lemma~\ref{lematec} one can assume without loss of generality that
$ \textrm{det}(A)=0$. In this situation, since $bc-ad=0$ we obtain, for $\lambda\neq 0$ and $\alpha:=\frac{a\beta}{\lambda}$ and all $p\geq 2$,
\begin{align*}
f_p=\frac{\lambda^p\prod_{i=0}^{p-1}\left(\frac{m(\mu_1+i)}{m(\mu_1+i-1)}-\frac{m(\mu_1)}{m(\mu_1-1)}+\alpha\right)}{\prod_{i=1}^p\left(\frac{m(\mu_1+i)}{m(\mu_1+i-1)}-\frac{m(\mu_1)}{m(\mu_1-1)}\right)\prod_{i=0}^{p-1}\left(\frac{m(\mu_1+i)}{m(\mu_1+i-1)}-\frac{m(\mu_1)}{m(\mu_1-1)}+\beta\right)},
\end{align*}
\begin{align*}
g_p=\frac{c\lambda^{p-1}\prod_{i=1}^{p-1}\left(\frac{m(\mu_1+i)}{m(\mu_1+i-1)}-\frac{m(\mu_1)}{m(\mu_1-1)}+\alpha\right)}{\prod_{i=1}^{p-1}\left(\frac{m(\mu_1+i)}{m(\mu_1+i-1)}-\frac{m(\mu_1)}{m(\mu_1-1)}\right)\prod_{i=1}^p\left(\frac{m(\mu_1+i)}{m(\mu_1+i-1)}-\frac{m(\mu_1)}{m(\mu_1-1)}+\beta\right)}.
\end{align*}
In the case that $\lambda=0$, for all $p\geq 2$, one has
\begin{align*}
f_p=\frac{a^p\beta^{p-1}}{\prod_{i=1}^p\left(\frac{m(\mu_1+i)}{m(\mu_1+i-1)}-\frac{m(\mu_1)}{m(\mu_1-1)}\right)\prod_{i=1}^{p-1}\left(\frac{m(\mu_1+i)}{m(\mu_1+i-1)}-\frac{m(\mu_1)}{m(\mu_1-1)}+\beta\right)},
\end{align*}
\begin{align*}
g_p=\frac{ca^{p-1}\beta^{p-1}}{\prod_{i=1}^{p-1}\left(\frac{m(\mu_1+i)}{m(\mu_1+i-1)}-\frac{m(\mu_1)}{m(\mu_1-1)}\right)\prod_{i=1}^p\left(\frac{m(\mu_1+i)}{m(\mu_1+i-1)}-\frac{m(\mu_1)}{m(\mu_1-1)}+\beta\right)}.
\end{align*}
The previous expressions can be obtained from the recursion formula. The computation corresponding to $\mu = \mu_2$ follows the same steps, leading to a second family of coefficients $s_{p,2}:=\left(\begin{matrix}
f_{p,2}
\\
g_{p,2}
\end{matrix}   \right) 
$, so we omit the details.

In conclusion, we have obtained two formal solutions of (\ref{eqhypsys}) in the form
$$y_1=\sum_{p\ge0}s_{p}z^{\mu_1+p},\qquad y_2=\sum_{p\ge0}s_{p,2}z^{\mu_2+p}.$$
The discussion on the convergence of such formal solutions can be made as in the general situation.

%%%%%%%%%%%%%%%%%%%%%%%%%%%%%%%%%%%%%%%%%%%%%%%555

\subsection{Planar systems. $B$ in Jordan form}\label{secnodiag}

We consider  the equation (\ref{eqhypsys}) with
$$
A:=\left(
\begin{matrix}
a&b
\\
c&d
\end{matrix}
 \right)\quad \textrm{and}	\quad  B:=\left(
\begin{matrix}
\mu_1&1
\\
0&\mu_1
\end{matrix}
 \right),
$$
where $\mu_1=\frac{m(\mu)}{m(\mu-1)}$ for some $\mu\in\C$ with $\hbox{Re}(\mu)\ge1$. We assume $\mu_1\neq\frac{m(\mu+p)}{m(\mu-1+p)}$ for any positive integer, so hypothesis (H1) holds for $\mu_1$. We proceed to search for a solution in the form of a formal Floquet power series, say (\ref{e416}).

Let us denote the coefficients of such solution by
$$
s_p:=\left(\begin{matrix}
f_p
\\
g_p
\end{matrix}   \right).
$$
The recursion for $s_p$ in \eqref{solution hypergeometric system} gives us, for all $p\geq 1$,
\begin{align}\label{particular solution Jordan case}
&\left(\frac{m(p+\mu)}{m(p+\mu-1)}-\mu_1\right)f_p=\left(a+\frac{c}{\frac{m(p+\mu)}{m(p+\mu-1)}-\mu_1}\right)f_{p-1}+\left(b+\frac{d}{\frac{m(p+\mu)}{m(p+\mu-1)}-\mu_1}\right)g_{p-1},
\\
&\left(\frac{m(p+\mu)}{m(p+\mu-1)}-\mu_1\right)g_p=cf_{p-1}+dg_{p-1}.\nonumber
\end{align}
Therefore, for  $\lambda=a+d$, since  $\mu_1=\frac{m(\mu)}{m(\mu-1)}$, we have for all $p\geq 1$,
\begin{align*}
&\left(\frac{m(p+1+\mu)}{m(p+\mu)}-\mu_1\right)\left(\frac{m(p+\mu)}{m(p+\mu-1)}-\mu_1\right)f_{p+1}
\\&
=
\left(a+\frac{c}{\frac{m(p+\mu+1)}{m(p+\mu)}-\mu_1}\right)\left(\frac{m(p+\mu)}{m(p+\mu-1)}-\mu_1\right)f_{p}+\left(b+\frac{d}{\frac{m(p+\mu+1)}{m(p+\mu)}-\mu_1}\right)\left(\frac{m(p+\mu)}{m(p+\mu-1)}-\mu_1\right)g_{p}
\\&=
\left(a+\frac{c}{\frac{m(p+\mu+1)}{m(p+\mu)}-\mu_1}\right)\left(\frac{m(p+\mu)}{m(p+\mu-1)}-\mu_1\right)f_{p}+\left(b+\frac{d}{\frac{m(p+\mu+1)}{m(p+\mu)}-\mu_1}\right)(cf_{p-1}+dg_{p-1})
\\&=
\left(a+\frac{c}{\frac{m(p+\mu+1)}{m(p+\mu)}-\mu_1}\right)\left(\frac{m(p+\mu)}{m(p+\mu-1)}-\mu_1\right)f_{p}+\left(b+\frac{d}{\frac{m(p+\mu+1)}{m(p+\mu)}-\mu_1}\right)cf_{p-1}
\\&\qquad+d\left(\left(\frac{m(p+\mu)}{m(p+\mu-1)}-\mu_1\right)f_p-\left(a+\frac{c}{\frac{m(p+\mu)}{m(p+\mu-1)}-\mu_1}\right)f_{p-1}\right)
\\&=
\left(\left(\lambda+\frac{c}{\frac{m(p+\mu+1)}{m(p+\mu)}-\frac{m(\mu)}{m(\mu-1)}}\right)\left(\frac{m(p+\mu)}{m(p+\mu-1)}-\frac{m(\mu)}{m(\mu-1)}\right)\right)f_{p}+(bc-da)f_{p-1}.
\end{align*}

In the case that $A$ and $B$ do not commute, the previous recursion determines the values of the coefficients of one Floquet solution. Otherwise, one can assume that $bc-ad=0$ (i.e. $\textrm{det}(A)=0$), in virtue of Lemma~\ref{lematec}, to arrive at  
\begin{align*}
g_0=0,\quad g_1=\frac{c}{\frac{m(\mu+1)}{m(\mu)}-\frac{m(\mu)}{m(\mu-1)}},\quad f_0=1\quad \textrm{and}\quad  f_1=\frac{a}{\frac{m(\mu+1)}{m(\mu)}-\frac{m(\mu)}{m(\mu-1)}}
\end{align*}
we obtain, for all $p\geq 2$,
\begin{align*}
f_p=\frac{\prod_{i=2}^{p}\left(\lambda+\frac{c}{\frac{m(\mu+i)}{m(\mu+i-1)}-\frac{m(\mu)}{m(\mu-1)}}\right)}{\prod_{i=1}^p\left(\frac{m(\mu+i)}{m(\mu+i-1)}-\frac{m(\mu)}{m(\mu-1)}\right)},
\end{align*}
\begin{align*}
g_p=\frac{c\left(\sum_{l=0}^{p-3}d^l\prod_{i=2}^{p-1-l}\left(\lambda+\frac{c}{\frac{m(\mu+i)}{m(\mu+i-1)}-\frac{m(\mu)}{m(\mu-1)}}\right)+d^{p-2}\lambda\right)}{\prod_{i=1}^p\left(\frac{m(\mu+i)}{m(\mu+i-1)}-\frac{m(\mu)}{m(\mu-1)}\right)},
\end{align*}
not appearing the series in the numerator for $p=2$. Both expressions can be verified by induction.

We have obtained a formal solution of (\ref{eqhypsys}) in the form
$$y_1=\sum_{p\ge0}s_pz^{\mu_1+p}.$$
The discussion on the convergence of such formal solution can be made as in the generic framework.

Finding a second solution of (\ref{eqhypsys}) is a much more delicate issue, due to the generality of moment derivation. In principal, if one is able to find a holomorphic function $H_m$, defined on some adequate domain $U\subseteq\C$ with $0\in\overline{U}$, which satisfies that 
\begin{equation}\label{e1141}
z\partial_m(H_m)(z)=\mu_1H_m(z)+z^{\mu}, \quad z\in U,
\end{equation}
then a second solution to the problem could be constructed, as we prove afterwards. However, many questions emerge from the existence of a function $H_m$ satisfying (\ref{e1141}). First, the values of the operator $\partial_m$ should be extended to a set of functions whose image is not only contained in $\C[[z]]$ (or $z^{\nu}\C[[z]]$ for some $\nu\in\C$ with $\hbox{Re}(\nu)\ge1$),  but on some broader set. As a matter of fact, such domain of definition would enlarge the set of analytic functions where it is currently well-defined: the set of analytic functions near the origin (by identifying the function with its Taylor expansion at the origin), the set of holomorphic functions defined on some sector with vertex at the origin which are the generalized sum of a formal power series (see~\cite{lastramichaliksuwinska21}), and the set of holomorphic functions defined on some sector with vertex at the origin which are the generalized multisum of a formal power series (see~\cite{lastramichaliksuwinska23}). In all the previous sets, $\partial_m$ has already been proved to be well-defined, with values in $\C[[z]]$. In the next section, we describe in a more general framework the existence of such function $H_m$ in some concrete situation, with importance in applications.

Let us assume the existence of an analytic function $H_m$ as described above, which satisfies (\ref{e1141}), at least from a symbolic point of view in two different situations. 

First, let us assume that $A\equiv 0$.

We first check that the symbolic expression
$$y_2(z)=\begin{pmatrix} H_m(z) \\ z^{\mu}\end{pmatrix}$$
is a solution of (\ref{eqhypsys}). Indeed, observe that
$$z\partial_my_2(z)=\begin{pmatrix} z\partial_m (H_m(z)) \\ \mu_1 z^{\mu}\end{pmatrix}.$$
On the other hand,
$$By_2(z)=\begin{pmatrix} \mu_1 & 1\\ 0 & \mu_1\end{pmatrix}\begin{pmatrix} H_m(z) \\ z^{\mu}\end{pmatrix}=\begin{pmatrix} \mu_1 H_m(z)+z^{\mu}\\ \mu_1 z^{\mu} \end{pmatrix},$$
which coincide in virtue of (\ref{e1141}). We conclude that 
$$\begin{pmatrix} 
z^{\mu} & H_m(z)\\
0 & 1
\end{pmatrix}
$$
is a solution to (\ref{eqhypsys}) with $A\equiv 0$. 

Let us conclude with the situation in which $A\nequiv0$. Let us consider the solution 
\begin{align*}
y_2(z)=\sum_{p=0}^{\infty}s_{p,1}z^{p+\mu}+\widetilde{H}_m(z)\sum_{p=0}^{\infty}s_{p,2}z^{p+\mu},
\end{align*}
 where we assume that $\widetilde{H}_m(z)$ is an holomorphic function satisfying
 \begin{align*}
 z\partial_m \left(\widetilde{H}_m(z)\sum_{p=0}^{\infty}s_{p,2}z^{p+\mu}\right)=
 \widetilde{H}_m(z)\sum_{p=0}^{\infty}s_{p,2}\frac{m(p+\mu)}{m(p+\mu-1)}z^{p+\mu}+\sum_{p=0}^{\infty}s_{p,2}z^{p+\mu}.
 \end{align*}
together with
 \begin{align*}
 s_{0,1}=\left(\begin{array}{c}
 0
 \\
 1
 \end{array}\right),\quad  s_{0,2}=\left(\begin{array}{c}
 1
 \\
 0
 \end{array}\right),
 \end{align*}
 and, for  $p\geq 1$,
 \begin{align}\label{eq:coeficientesAneq0}
 s_{p,2}=\left(\frac{m(p+\mu)}{m(p+\mu-1)}I-B\right)^{-1}As_{p-1,2},
  \quad
 \textrm{and}\quad
  s_{p,1}=\left(\frac{m(p+\mu)}{m(p+\mu-1)}I-B\right)^{-1}\left(As_{p-1,1}-s_{p,2}\right).
 \end{align}
 Indeed, observe that
 \begin{align}\label{eq:derivada ejemplo}
 z\partial_m y_2(z)
 &
 =\sum_{p=0}^{\infty}s_{p,1}\frac{m(p+\mu)}{m(p+\mu-1)}z^{p+\mu}
 + \widetilde{H}_m(z)\sum_{p=0}^{\infty}s_{p,2}\frac{m(p+\mu)}{m(p+\mu-1)}z^{p+\mu}+\sum_{p=0}^{\infty}s_{p,2}z^{p+\mu}
 \\&\nonumber
 =\sum_{p=1}^{\infty}s_{p,1}\frac{m(p+\mu)}{m(p+\mu-1)}z^{p+\mu}
 + \widetilde{H}_m(z)\sum_{p=1}^{\infty}s_{p,2}\frac{m(p+\mu)}{m(p+\mu-1)}z^{p+\mu}+\sum_{p=1}^{\infty}s_{p,2}z^{p+\mu}
 \\&\nonumber\quad
 +\left(
s_{0,1}\frac{m(\mu)}{m(\mu-1)}
 + \left(\widetilde{H}_m(z)\frac{m(\mu)}{m(\mu-1)}+1\right)s_{0,2}\right)z^{\mu}.
 \end{align}
 On the other hand, by \eqref{eq:coeficientesAneq0}
 \begin{align}\label{eq:matrices ejemplo}
 (zA+B)y_2(z)&=
 \sum_{p=0}^{\infty}As_{p,1}z^{p+\mu+1}+\widetilde{H}_m(z)\sum_{p=0}^{\infty}As_{p,2}z^{p+\mu+1}+\sum_{p=0}^{\infty}Bs_{p,1}z^{p+\mu}+\widetilde{H}_m(z)\sum_{p=0}^{\infty}Bs_{p,2}z^{p+\mu}
  \\&\nonumber
 =
 \sum_{p=1}^{\infty}As_{p-1,1}z^{p+\mu}+\widetilde{H}_m(z)\sum_{p=1}^{\infty}As_{p-1,2}z^{p+\mu}+\sum_{p=1}^{\infty}Bs_{p,1}z^{p+\mu}+\widetilde{H}_m(z)\sum_{p=1}^{\infty}Bs_{p,2}z^{p+\mu}
 \\&\nonumber
 \quad + \left(Bs_{0,1}+\widetilde{H}_m(z)Bs_{0,2}\right)z^{\mu}
  \\&\nonumber
 =
 \sum_{p=1}^{\infty}\left(\left(\frac{m(p+\mu)}{m(p+\mu-1)}I-B\right)s_{p,1}+s_{p,2}\right)z^{p+\mu}+\widetilde{H}_m(z)\sum_{p=1}^{\infty}\left(\frac{m(p+\mu)}{m(p+\mu-1)}I-B\right)s_{p,2}z^{p+\mu}
 \\&\nonumber
 \quad +\sum_{p=1}^{\infty}Bs_{p,1}z^{p+\mu}+\widetilde{H}_m(z)\sum_{p=1}^{\infty}Bs_{p,2}z^{p+\mu} + \left(Bs_{0,1}+\widetilde{H}_m(z)Bs_{0,2}\right)z^{\mu}
  \\&\nonumber
 =
 \sum_{p=1}^{\infty}\left(\left(\frac{m(p+\mu)}{m(p+\mu-1)}I\right)s_{p,1}+s_{p,2}\right)z^{p+\mu}+\widetilde{H}_m(z)\sum_{p=1}^{\infty}\left(\frac{m(p+\mu)}{m(p+\mu-1)}I\right)s_{p,2}z^{p+\mu}
 \\&\nonumber
 \quad + \left(Bs_{0,1}+\widetilde{H}_m(z)Bs_{0,2}\right)z^{\mu}.
 \end{align}
 Next, note that
 \begin{align*}
 Bs_{0,1}=\mu_1 s_{0,1}+s_{0,2}=\frac{m(\mu)}{m(\mu-1)}  s_{0,1}+s_{0,2}
 \quad \textrm{and}
 \quad Bs_{0,2}=\mu_1 s_{0,2}=\frac{m(\mu)}{m(\mu-1)} s_{0,2}.
 \end{align*}
  Combining this with \eqref{eq:derivada ejemplo} and  \eqref{eq:matrices ejemplo}, gives us
 \begin{align*}
  (zA+B)y_2(z)&=
   \sum_{p=1}^{\infty}\left(\left(\frac{m(p+\mu)}{m(p+\mu-1)}I\right)s_{p,1}+s_{p,2}\right)z^{p+\mu}+\widetilde{H}_m(z)\sum_{p=1}^{\infty}\left(\frac{m(p+\mu)}{m(p+\mu-1)}I\right)s_{p,2}z^{p+\mu}
 \\&\nonumber
 \quad + \left(\frac{m(\mu)}{m(\mu-1)}  s_{0,1}+\left(1+\widetilde{H}_m(z)\frac{m(\mu)}{m(\mu-1)} \right)s_{0,2}\right)z^{\mu}
 \\
 &=
   \sum_{p=1}^{\infty}\frac{m(p+\mu)}{m(p+\mu-1)}s_{p,1}z^{p+\mu}+\sum_{p=1}^{\infty}s_{p,2}z^{p+\mu}+\widetilde{H}_m(z)\sum_{p=1}^{\infty}\frac{m(p+\mu)}{m(p+\mu-1)}s_{p,2}z^{p+\mu}
 \\&\nonumber
 \quad + \left(\frac{m(\mu)}{m(\mu-1)}  s_{0,1}+\left(1+\widetilde{H}_m(z)\frac{m(\mu)}{m(\mu-1)} \right)s_{0,2}\right)z^{\mu}
 \\&
 =  z\partial_m y_2(z).
 \end{align*}

\begin{remark} Observe that $H_m$ in the case of $A\equiv0$ coincides with $\widetilde{H}_m(z) z^{\mu}$ for $A\not\equiv0$.
\end{remark}

\subsection{Functions of a matrix and confluent hypergeometric systems}

In Section~\ref{secfunmat}, a definition generalizing $z^B$, for any diagonalizable matrix $B\in\C^{n\times n}$ under the assumptions of Theorem~\ref{teopral1} with $k=n$, has been put forward. The columns of this matrix, $z_m^B$, provide $n$ solutions to the system
$$z\partial_my=By,$$
provided that a function $H_m$ as in (\ref{e1141}) exists. For this purpose, we make the following assumption:

\begin{itemize}
\item[(H3)] For every $\mu\in\C^{\star}$, with $\hbox{Re}(\mu)\ge1$, there exists a sequence of functions $(H_{m,p,\mu})_{p\ge0}$ satisfying the following recursive family of moment differential equations on some domain $U\subseteq\C$, with $0\in\overline{U}$. Define $H_{m,1,\mu}:=z^{\mu}$ and for all $p\ge1$
$$z\partial_m(H_{m,p+1,\mu})(z)=\mu_1H_{m,p+1,\mu}(z)+H_{m,p,\mu}(z), \quad z\in U,$$
where $\mu_1=\frac{m(\mu)}{m(\mu-1)}$.
\end{itemize}

Taking into account the facts observed in the previous section, the last main result of the present work provides the symbolic solution of 
$$z\partial_m=By$$
for any $B\in\C^{n\times n}$, under Hypothesis (H3).  

\begin{theorem}\label{teopral2}
Let $B\in\C^{n\times n}$ whose Jordan canonical form is given by $J=\hbox{diag}(J_1,\ldots,J_s)$, where $J_k\in\C^{n_k\times n_k}$ for $1\le k\le s$ corresponds to a Jordan block 
$$J_k=\begin{pmatrix}
\mu_k & 1       & 0      & \cdots  & 0 \\
0       & \mu_k & 1      & \cdots  & 0 \\
\vdots  & \vdots  & \vdots & \ddots  & \vdots \\
0       & 0       & 0      & \mu_k & 1      \\
0       & 0       & 0      & 0       & \mu_k
\end{pmatrix},$$
for some positive integer $n_k$ and $\mu_k\in\C$ such that there exists $\tilde{\mu}_k\in\C$ with $\hbox{Re}(\tilde{\mu}_k)\ge1$ and $\mu_k=\frac{m(\tilde{\mu}_k)}{m(\tilde{\mu}_k-1)}$.  We assume that for all $1\le k_1,k_2\le s$ with $k_1\neq k_2$ and all positive integer $n$ one has that 
$$\mu_{k_1}\neq\frac{m(\tilde{\mu}_{k_2}+n)}{m(\tilde{\mu}_{k_2}+n-1)},$$

Then, there exist $s$ formal solutions solutions of 
\begin{equation}\label{e1176}
z\partial_my=By
\end{equation}
in the form of a Floquet formal power series, say $C_{1},\ldots,C_s$. Moreover, if (H3) holds (substituting $\mu$ and $\mu_1$ in (\ref{e1141}) by $\tilde{\mu}_k$ and $\mu_k$, respectively) for all $1\le k\le s$, then there exist $n-s$ symbolic solutions: $C_{1,1},\ldots,C_{1,n_1-1}$ associated to $C_1$, $\ldots$, and $C_{s,1},\ldots, C_{s,n_s-1}$ associated to $C_s$ ($n_1+\ldots+n_s=n$). The symbolic matrix with columns given by the previous vectors
$$z_m^{B}:=(C_{1}(z),C_{1,1}(z),\ldots,C_{1,n-1}(z),\ldots,C_{s}(z),C_{s,1}(z),\ldots, C_{s,n_s-1}(z))$$ 
is a symbolic solution to (\ref{e1176}).
\end{theorem}

\begin{remark} The explicit form of the columns in $z_m^B$ is provided in the proof of Theorem~\ref{teopral2}. \end{remark}

\begin{proof}
First, observe from Proposition~\ref{prop1} that one can reduce the problem to the case in which $B$ is already given by its Jordan associated matrix $J$. We will construct the solution when there is only one block in $J$, i.e. $s=1$, without loss of generality. The procedure for $s\ge2$ follows an analogous reasoning, completing the solutions according to the classical procedure in the framework of systems of differential equations, namely adding zeros at the coordinates which do not correspond to the block under consideration.

After taking into account the previous considerations, we rewrite the problem in the following form for the sake of a clear reading:
\begin{equation}\label{e1202}
z\partial_my=Jy,
\end{equation}
with $J\in\C^{n\times n}$, and where 
$$J=\begin{pmatrix}
\mu_1 & 1       & 0      & \cdots  & 0 \\
0       & \mu_1 & 1      & \cdots  & 0 \\
\vdots  & \vdots  & \vdots & \ddots  & \vdots \\
0       & 0       & 0      & \mu_1 & 1      \\
0       & 0       & 0      & 0       & \mu_1
\end{pmatrix},$$
for some $\mu_1\in\C$ such that $\frac{m(\mu)}{m(\mu-1)}=m(\mu_1)$ for certain $\mu\in\C$ with $\hbox{Re}(\mu)\ge1$. A first solution of (\ref{e1202}) is constructed by means of the eigenvetor $s=(1,0,\ldots,0)$ associated to the eigenvalue $\mu_1$ of $J$ (see Definition~\ref{defin1}) providing the first column of $z_m^J$ determined by $C_1=(z^{\mu},0,\ldots,0)^{T}$. 

Assumption (H3) guarantees the existence of a domain $U\subseteq\C$ with $0\in\overline{U}$, and functions $H_1(z)=z^{\mu},H_2\ldots,H_{n}$ such that
$$z\partial_m(H_{p+1})(z)=\mu_1H_{p+1}(z)+H_{p}(z), \quad z\in U,$$
for all $1\le p\le n-1$. The remaining $n-1$ solutions determine $C_2,\ldots,C_n$ defined by
\begin{multline*}
C_2=(H_2(z),z^{\mu},0,\ldots,0)^T\\
C_3=(H_3(z),H_2(z),z^{\mu},0,\ldots,0)^T\\
\cdots\\
C_{n}=(H_{n}(z),H_{n-1}(z),\ldots,H_2(z),z^{\mu})^T.
\end{multline*}
It is straight to check that the matrix with columns $C_1,\ldots,C_n$ is a solution of (\ref{e1202}).
\end{proof}

\begin{remark} Observe from the construction of $z_m^B$ in the previous Theorem generalizes that of Definition~\ref{defin1}, associated to a diagonalizable matrix $B$.
\end{remark}

\begin{remark} The hypothesis $\hbox{Re}(\tilde{\mu}_k)\ge1$ can be obviated from the hypotheses of Theorem~\ref{teopral2} adapting the hypotheses and following analogous arguments as those in Section~\ref{sec35}.
\end{remark}

We conclude the present work by describing two concrete situations in which hypothesis (H3) holds for some holomorphic sequence of holomorphic functions $(H_{m,p,\mu})_{p\ge0}$.

\begin{example}
Let us consider the moment sequence $m=(p!)_{p\ge0}$ which can be extended to the complex numbers of positive real part by $m(z)=\Gamma(1+z)$. The theory described is then reduced to the classical theory of differential equations, due to $\partial_m$ turns out to be the classical derivative. The functions $H_p:=H_{m,p,\mu}$ are given by
$$H_{p}(z)=\frac{1}{p!}z^{\mu}\log^{p}(z),\quad p\ge 0.$$
These functions are holomorphic on $\C\setminus(-\infty,0]$. Observe that $H_1(z)=z^{\mu}$ and for every $p\ge1$ one has that
\begin{multline*}
z(H_{p+1}(z))'=\frac{1}{(p+1)!}z(z^{\mu}\log^{p+1}(z))'=\frac{\mu}{(p+1)!}z^{\mu}\log^{p+1}(z)+\frac{1}{p!}z^{\mu}\log^p(z)\\
=\mu H_{p+1}+H_p.
\end{multline*}
Recall that 
$$\mu_1=\frac{m(\mu)}{m(\mu-1)}=\frac{\Gamma(1+\mu)}{\Gamma(\mu)}=\mu.$$
\end{example}

\begin{example}
Assume $\mu$ is a positive integer. Let us fix $m=([p]_q^{!})_{p\ge0}$, for some $q>1$. We recall from Section~\ref{secapap} that this sequence is determined by the evaluation of $q$-Gamma function $z\mapsto\Gamma_q(1+z)$ at the nonnegative integers. In addition to this, $\partial_m$ corresponds to Jackson's $q$-derivative $D_q$. The existence of the sequence of functions $(H_{m,p,\mu})_{p\ge0}$ is guaranteed. Let us denote it $(H_p)_{p\ge0}$ for simplicity. 

Indeed, $H_1(z)=z^{\mu}$ which is an entire function. Moreover, the recursive family of $q-$difference equations can be written in this context in the form
\begin{equation}\label{e1253}
H_{p+1}(qz)=(1+(q-1)[\mu]_q)H_{p+1}(z)+(q-1)H_p,\quad p\ge1.
\end{equation}

At this point, we apply the following results from~\cite{lastraremy}, adapted to our needs.

\begin{lemma}[Proposition 2.2(2),~\cite{lastraremy} and~\cite{sau}]\label{lemafin1}
Let $c\in\C\setminus\{0,1\}$. The equation 
$$y(qz)=cy(z)$$
admits 
$$y(z)=\frac{\Theta_q(z)}{\Theta_{q}(z/c)}$$
as a meromorphic solution on $\C^{\star}$, with simple zeros and poles located at $\{q^{n}:n\in\mathbb{Z}\}$ and $\{cq^{n}:n\in\mathbb{Z}\}$, respectively.
\end{lemma}
In the previous result, $\Theta_q(\cdot)$ stands for Jacobi Theta function.

\begin{lemma}[Proposition 2.15,~\cite{lastraremy}]\label{lemafin2}
Let $R(z)$ be a meromorphic function on $\C^{\star}$ whose poles are positive real numbers. Then, the functional equation
$$y(qz)=y(z)+R(z)$$
admits a meromorphic solution on $\C^{\star}$ whose poles are positive real numbers
\end{lemma}

The proof of the previous result also provides a constructive solution to the $q$-difference equation under study.

Finally, we see from Section 2.3~\cite{lastraremy} that if $M(z)$ is a meromorphic function in $\C^{\star}$ which solves the equation $y(qz)=cy(z)$, and $r(z)$ is any given meromorphic function in $\C^{\star}$, then the meromorphic solution on $\C^{\star}$ of the equation $y(qz)=y(z)+M^{-1}(qz)r(z)$ is a meromorphic solution on $\C^{\star}$ of the equation
$$y(qz)=cy(z)+r(z).$$

Let $c=1+(q-1)[\mu]_q)$ in Lemma~\ref{lemafin1}, and let $R(z)=H_p(z)$ in Lemma~\ref{lemafin2}. The application of Lemma~\ref{lemafin1} and Lemma~\ref{lemafin2} allow us to conclude the existence of holomorphic functions $H_p$ defined on $\C\setminus[0,\infty)$ which satisfy hypothesis (H3).
\end{example}

%%%%%%%%%%%%%%%%%%%%%%%%%%%%%%%%%%%%%%%%%%%%%%%%%%%%%%%%%%%%%%%%%%%%%%%%%5

\end{document}